\documentclass{article}

\usepackage{amssymb}
\usepackage{latexsym}
\usepackage{amsthm}
\usepackage{epsfig}
\usepackage{float}
\usepackage{amsmath}
\usepackage{ifthen}
\usepackage{xspace}
\usepackage[nobysame,alphabetic]{amsrefs}
\usepackage{srcltx}
\usepackage{amscd}

\def\a{\alpha}

\def\b{\beta}
\def\C{\mathcal{C}}
\def\g{\gamma}
\def\G{\Gamma}

\def\H{\mathcal{H}}

\def\l{\lambda}

\def\O{\mathcal{O}}

\def\R{\mathbb{R}}
\def\s{\sigma}
\def\S{\Sigma}
\def\T{\mathcal{T}}
\def\Z{\mathbb{Z}}
\def\d{\partial}

\def\union{\cup}

\def\e{\epsilon}

\def\P{\mathbb{P}}
\def\w{\omega}

\def\mun{\mu^{(n)}}

\def\PML{\mathcal{PML}}
\def\PMF{\mathcal{PMF}}
\def\fmin{\mathcal{F}_{min}}
\def\lmin{\mathcal{L}_{min}}
\def\UE{\mathcal{UE}}
\def\Ghat{\widehat \Gamma}

\def\Grel{\widehat G}
\def\gz{G^{\Z_+}}

\def\Kcc{Q}
\def\ccc{q}
\def\Kconj{K}
\newcommand{\K}[1]{K_{#1}}

\def\fix{\text{fix}}

\def\axis{\a}
\def\axisa{\a}
\def\axisb{\b}
\def\Thin{\text{Thin}}
\def\qcenter{\text{Center}}
\def\pA{pseudo-Anosov\xspace}
\def\mcg{mapping class group\xspace}

\def\le{\leqslant}
\def\ge{\geqslant}

\newcommand{\norm}[1]{|#1|}
\newcommand{\nhat}[1]{|\widehat{#1}|}
\newcommand{\dhat}[1]{\widehat d (#1)}
\newcommand{\dc}[1]{d_\C(#1)}
\newcommand{\cfix}[1]{\fix_{6\delta}(#1)}
\newcommand{\horo}[2]{\mathcal{O}_#2(#1)}
\newcommand{\geo}[1]{[1,#1]}
\renewcommand{\path}[2]{[#1,#2]}

\newtheorem{theorem}{Theorem}[section]
\newtheorem{lemma}[theorem]{Lemma}

\newtheorem{proposition}[theorem]{Proposition}

\newtheorem*{theorem:horoball}{Theorem \ref{theorem:horoball}}

\theoremstyle{definition}

\newcounter{case}

\newtheorem*{theorem:conjugacy}{Theorem \ref{theorem:conjugacy}}
\newtheorem*{lemma:small}{Lemma \ref{lemma:small}}

\restylefloat{figure}

\begin{document}


\title{Random walks on the mapping class group}
\author{Joseph Maher\footnote{email: joseph.maher@csi.cuny.edu}}
\date{\today}

\maketitle

\begin{abstract}
We show that a random walk on the mapping class group of an orientable
surface gives rise to a pseudo-Anosov element with asymptotic
probability one. Our methods apply to many subgroups of the mapping
class group, including the Torelli group.


Subject code: 
37E30, 
20H10, 
60G50, 
20F65. 

\end{abstract}

\tableofcontents

\section{Introduction}

Let $\S$ be an orientable surface of finite type, i.e. a surface of
genus $g$ with $p$ marked points, usually referred to as punctures.
The mapping class group of $\S$ consists of orientation preserving
diffeomorphisms which preserve the punctures, modulo those isotopic to
the identity. Thurston \cite{thurston2} showed that all elements of
the mapping class group are periodic, reducible or pseudo-Anosov. In
this paper we show that a random walk on the mapping class group gives
rise to a pseudo-Anosov element with asymptotic probability one. In
fact we obtain a more general result which we now describe.

Let $\mu$ be a probability distribution on the mapping class group
$G$.  A random walk on $G$ is a Markov chain on $G$ with transition
probabilities given by left translation of $\mu$, i.e. the probability
that you go from $x$ at time $n$, to $y$ at time $n+1$, is $p(x,y) =
\mu(x^{-1}y)$, and we shall assume that we start at the identity at
time zero. The path space for the random walk is the probability space
$(\gz, \P)$, where the product $\gz$ is the collection of all sample
paths, and the measure $\P$ is determined by $\mu$. If $w$ is a sample
path, then we will write $w_n$ for the location of the sample path at
time $n$, and the distribution of $w_n$ is given by the $n$-fold
convolution of $\mu$.  The support of the random walk is the
semi-group generated by the support of $\mu$. We shall always assume
that the group generated by the support of $\mu$ is a non-elementary
subgroup of the mapping class group, i.e. it contains a pair of \pA
elements with distinct fixed points in $\PML$. The mapping class group
is finitely generated, so a choice of (finite) generating set gives
rise to a word metric on the group, where the length of an element is
the shortest length of any word in the generators equal to the group
element, and two different choices of finite generating sets give
quasi-isometric word metrics.  We may also consider the word metric on
the mapping class group with respect to an infinite generating set,
though in this case the resulting metric need not be proper, and is
often referred to as a relative metric.  Masur and Minsky \cite{mm1}
have shown that there is a choice of infinite generating set,
consisting of a finite generating set union a particular collection of
subgroups, which gives rise to a metric which is $\delta$-hyperbolic.
In this case, we say that the mapping class group is weakly relatively
hyperbolic. In fact, the mapping class group with this particular
relative metric is quasi-isometric to the complex of curves.  This
gives two different ways to measure lengths of elements of the mapping
class group, as we can measure their length using word length with
respect to a finite generating set, or we can measure their length in
the relative metric.  The relative length of a group element is its
distance from the identity in the relative metric. We show that $w_n$
is conjugate to an element of bounded relative length with asymptotic
probability zero, i.e. the probability that $w_n$ is conjugate to an
element of bounded relative length tends to zero as $n$ tends to
infinity.

\begin{theorem} 
\label{theorem:rw}
Consider a random walk on the mapping class group of an orientable
surface of finite type, which is not a sphere with three or fewer
punctures, and let $w_n$ be the location of the random walk after $n$
steps.  If the group generated by the support of the random walk is
non-elementary, then for any constant $B$, $w_n$ is conjugate to an
element of relative length at most $B$ with asymptotic probability
zero.
\end{theorem}

There is a constant $B$, depending on the surface, such that any
mapping class group element which is not \pA is conjugate to an
element of relative length at most $B$, as shown in Lemma
\ref{lemma:bounded length}, and so this implies that $w_n$ is
pseudo-Anosov with asymptotic probability one.  Rivin \cite{rivin} and
Kowalski \cite{kowalski} have also shown that the nearest neighbour
random walk on the mapping class group gives rise to a pseudo-Anosov
element with asymptotic probability one, as part of a broader
investigation of random walks on groups. Their methods apply whenever
the support of the random walk maps onto $Sp(2n,\Z)$, and have the
advantage that they may be used to construct explicit lower bounds for
the proportion of elements which are pseudo-Anosov at time $n$, and
furthermore, show that the set of \pA elements is transient for these
random walks.  However, our methods apply to random walks supported on
more general subgroups of the mapping class group. For example, the
Torelli group is the subgroup of the mapping class group which acts
trivially on the homology of the surface. This is a normal subgroup
whose limit set is dense in the Thurston boundary, and so is not
contained in a non-trivial centralizer.  Therefore the nearest
neighbour random walk on a Cayley graph for the Torelli group gives
rise to a pseudo-Anosov element with asymptotic probability one.

We now give an outline of the main argument. Let $R$ be a set of
elements in the mapping class group which are all conjugate to
elements of bounded relative length, for example all of the non-\pA
elements of the mapping class group. The basic idea is to consider the
distribution of elements of $R$ inside the mapping class group, and
show that random walks end up travelling through regions in which the
density of elements of $R$ tends to zero. The long time behaviour of
random walks in the mapping class group is described by harmonic
measure on $\PML$, the space of projective measured laminations. This
space may be thought of as a boundary for the mapping class group,
following work of Masur and Minsky \cite{mm1} and Klarreich
\cite{klarreich}, and we make this precise in Section
\ref{section:random}. By work of Kaimanovich and Masur \cite{km},
sample paths converge to the boundary with probability one, and this
gives rise to a measure on $\PML$, called harmonic measure. The
measure of a subset of the boundary is the probability that a sample
path will converge to a lamination lying in that set, and this measure
depends on the choice of probability distribution $\mu$ used to define
the random walk.  Given a subset $X$ of the mapping class group, we
can take the harmonic measure of its limit set $\overline X$ in
$\PML$. A sample path of the random walk is recurrent on $X$ if it
hits $X$ infinitely often. If a sample path is recurrent on $X$, and
it converges to the boundary, then it converges to a point in
$\overline X$. As sample paths converge to the boundary with
probability one, this means that the probability that a sample path is
recurrent on $X$ is a lower bound on the harmonic measure of
$\overline X$. In particular, if the harmonic measure of the limit set
is zero, then the random walk is transient on $X$, i.e. a sample path
hits $X$ only finitely many times with probability one.

In general the limit set of $R$ need not have harmonic measure zero,
for example the set of all elements of the mapping class group which
are not pseudo-Anosov has a limit set consisting of all of $\PML$, and
so has harmonic measure one.  However, we can consider $R_k$, the set
of all the elements in $R$ which are word length distance (not
relative distance) at most $k$ from some other element of $R$. We show
that the limit set of $R_k$ is contained in the union of the limit
sets of the centralizers of elements of the mapping class group of
length at most $k$. We then show that this set has harmonic measure
zero, assuming that the group $H$ generated by the support of the
random walk is non-elementary, and every non-trivial element has a
centralizer which has infinitely many images under $H$.  Each element
of $R \setminus R_k$ lives inside a ball of elements not in $R$ of
radius at least $k$, and so there is an upper bound for the
probability that a sample path hits an element of $R \setminus R_k$.
This bound tends to zero as $k$ tends to infinity, so the probability
that a sample path is not pseudo-Anosov tends to zero as the length of
the path tends to infinity. Finally, if there are centralizers with
finitely many images under the group $H$ generated by the
support of the random walk, we show that we can map $H$ to the mapping
class group of a surface covered by $\S$, whose mapping class group
does satisfy the condition on images of centralizers.

An important part of showing that the limit set of $R_k$ is a finite
union of centralizers, is showing that if an element $g$ is conjugate
to a relatively short element $s$, i.e. $g = wsw^{-1}$, and if $w$ is
chosen to be the shortest conjugating word, then the path $wsw^{-1}$
is actually quasi-geodesic, for quasi-geodesic constants that depend
on the relative length of $s$, but not of $g$.  This follows from
showing that the mapping class group has relative conjugacy bounds,
which means that the relative length of a shortest conjugating word is
bounded in terms of the relative lengths of the two conjugate words.
Using further work of Masur and Minsky \cite{mm2}, we show that the
mapping class group has relative conjugacy bounds by considering the
action of the mapping class group on the complex of curves. The
relative conjugacy bound property is equivalent to discreteness for
the action of the mapping class group on the complex of curves, in the
following sense: hyperbolic isometries have a minimal translation
distance, and elliptic isometries have centralizers which act coarsely
transitively on their coarse fixed sets.

The paper is structured as follows. In Section \ref{section:defs} we
review some well known definitions and results.
In Section \ref{section:conjugacy} we show that the mapping class
group has relative conjugacy bounds. In Section \ref{section:coarse}
we describe the distribution of $k$-dense reducible elements in the
mapping class group, and finally in Section \ref{section:random} we
prove the main result on random walks.

\subsection{Acknowledgements}

I would like to thank Nathan Dunfield and Howard Masur for useful
advice. I would especially like to thank the referees for pointing out
many errors in an earlier version, supplying a correct proof of Lemma
\ref{lemma:ref}, and for suggestions on how to re-organize the paper.
I would also like to thank Danny Calegari, Daniel Groves and Jason
Manning for helpful conversations. Part of this paper was written
while I was at the California Institute of Technology, and part while
I was supported by a Postdoctoral Fellowship from the Centre de
recherches math\'ematiques and the Institut des sciences
math\'ematiques in Montr\'eal.  This work was also partially supported
by NSF grant DMS-070674.

\section{Preliminary definitions} \label{section:defs}

In this section we review some of the definitions and results we will
use in the main part of the paper, and fix some notation.

Let $\S$ be the closed orientable surface of genus $g$ with $n$ marked
points, also known as punctures. We will write $\S_{g,n}$ if we need
to explicit refer to the genus $g$, and the number of punctures
$n$. The mapping class group of $\S$ is the group of orientation
preserving diffeomorphisms of the surface which preserve the set of
punctures, modulo those isotopic to the identity. We say a surface is
\emph{sporadic} if it is a sphere with at most four punctures, or a
torus with at most one puncture.  If $\S$ is a sphere with three or
fewer punctures, then the mapping class group of $\S$ is finite, and
does not contain any \pA elements.  If $\S$ is a torus with at most
one puncture, or a four punctured sphere, then the mapping class group
is commensurable with $SL(2,\Z)$, and random walks on $SL(2,\Z)$ are
well understood by work of Furstenberg \cite{furstenberg}.  So for the
remainder of this paper we will assume that the surface $\S$ is
\emph{not} sporadic. However, we will at times need to consider
subsurfaces of $\S$ which may be sporadic, so certain results will be
needed for surfaces which include a torus with one puncture, or a four
punctured sphere.

Dehn \cite{dehn} and Lickorish \cite{lickorish} showed that the
mapping class group is finitely generated. In fact, we can describe an
explicit generating set using Dehn twists. A Dehn twist is a map from
the surface to itself defined by cutting the surface along an
essential simple closed curve, and then gluing the two boundary
components back together using a full twist. We say a simple closed
curve in the surface is \emph{essential} if it does not bound a
subsurface which is either a disc, or a disc containing a single
puncture. A Dehn twist represents a non-trivial element of the mapping
class group. A Dehn twist is supported in a neighbourhood of the
simple closed curve used to define it, and takes a transverse arc in a
regular neighbourhood of the simple closed curve to one which winds
once around the annular regular neighbourhood. This is illustrated in
Figure \ref{picture31} below.

\begin{figure}[H]
\begin{center}
\epsfig{file=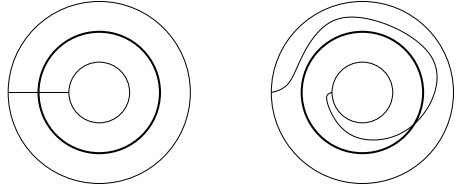, height=70pt}
\end{center}
\caption{A Dehn twists in a simple closed curve.} \label{picture31}
\end{figure}

The collection of Dehn twists in the simple closed curves shown below
in Figure \ref{picture21} generate the mapping class group of a closed
surface with no punctures. If the surface has punctures, then there is
a finite collection of Dehn twists which generate the finite index
subgroup of the mapping class group which fixes each puncture, and
this may be extended to a finite generating set for the whole group by
adding elements which permute the punctures. These elements may be
chosen to be half Dehn twist in simple closed curves that bound discs
containing exactly two punctures, see for example Birman
\cite{birman}*{Chapter 4}.

\begin{figure}[H]
\begin{center}
\epsfig{file=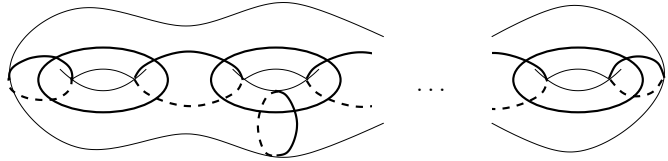, height=60pt}
\end{center}
\caption{Dehn twists which generate the mapping class group.} \label{picture21}
\end{figure}

We shall think of elements of the mapping class group as functions on
surfaces, so we will write $a b$ to denote $ab(x) = a(b(x))$ for all
$x \in \S$.

The collection of essential simple closed curves in the surface may be
made in to a simplicial complex, called the the \emph{complex of
  curves}, which we shall denote $\C(\S)$.  The vertices of this
complex are isotopy classes of simple closed curves in $\S$, and a
collection of vertices spans a simplex if representatives of the
curves can be realised disjointly in the surface.  The complex of
curves is a finite dimensional simplicial complex, but it is not
locally finite. We will write $\C_0(\S)$ to denote the vertices of the
simplicial complex $\C(\S)$, which is the set of isotopy classes of
simple closed curves.  We will write $\dc{x,y}$ for the distance in
the one-skeleton between two vertices $x$ and $y$ of the complex of
curves. We will always consider the complex of curves to have a
basepoint $x_0$, which we shall take to be one of the curves
corresponding to a standard generating set for the mapping class
group. The mapping class group acts by simplicial isometries on the
complex of curves. For certain sporadic surfaces the definition above
produces a collection of disconnected points, and so a slightly
different definition is used. If the surface is a torus with at most
one puncture, then two vertices are connected by an edge if the
corresponding simple closed curves may be isotoped to intersect
transversely exactly once. If the surfaces is a four punctured sphere,
then two vertices are connected by an edge if the corresponding simple
closed curves may be isotoped to intersect transversely in two points.
In these cases, the resulting curve complex is isomorphic to the Farey
graph.

A geodesic metric space is \emph{$\delta$-hyperbolic} if every
geodesic triangle is $\delta$-slim, i.e. each edge is contained in a
$\delta$-neighbourhood of the other two.  Masur and Minsky \cite{mm1}
have shown that the complex of curves is $\delta$-hyperbolic.

\begin{theorem}{\cite{mm1}*{Theorem 1.1}}
Let $\S$ be an oriented surface of finite type. The curve complex
$\C(\S)$ is $\delta$-hyperbolic, and has infinite diameter unless $\S$
is a sphere with three or fewer punctures.
\end{theorem}

An explicit bound for $\delta$ in terms of the complexity of the surface is
given by Bowditch \cite{bb}.

According to Thurston's classification of surface automorphisms, there
are three sorts of elements of the mapping class group. We say that an
element $g$ of the mapping class group is \emph{periodic} if it has
finite order. We say $g$ is \emph{reducible}, if it fixes a disjoint
collection of simple closed curves in $\S$. The periodic and reducible
elements act on the complex of curves as elliptic isometries. All
other elements are \emph{pseudo-Anosov}, and act on the complex of
curves as hyperbolic isometries.

Let $G$ be a finitely generated group, together with a symmetric
generating set $A$. For each element $g$ of $G$, we define the word
length of $g$ to be the length of the shortest word representing $g$
in the generating set $A$. This length function on $G$ induces a
left-invariant distance metric on $G$, called the \emph{word metric}.
This metric may also be obtained by forming the \emph{Cayley graph}
$\G$ for the group.  The Cayley graph is the graph whose vertices
consist of elements of $G$, with two group elements $a$ and $b$
connected by an edge if $a^{-1}b$ is a member of the generating set.
If we set the length of each edge equal to $1$ then the distance
between two vertices in the Cayley graph is the same as the word
metric distance.  Different choices of generating set give different,
but quasi-isometric, word metrics. We will assume we have fixed some
generating set for the mapping class group, and we will write
$\norm{x}$ for the word length of $x$, and $d(x,y)$ for the distance
from $x$ to $y$ in the word metric on $G$.  A geodesic in the word
metric on $G$ is a word of minimal length in the generating set.

Given a group $G$, and a collection of subgroups $\H = \{ H_i \}_{i
  \in I }$, we define the \emph{relative length} of a group element
$g$ to be the length of the shortest word in the typically infinite
generating set $A \cup \H$.  This defines a metric on $G$ called the
\emph{relative metric}, which depends on the choice of subgroups $\H$.
We will write $\Grel$ to denote the group $G$ with the relative
metric, which we shall also refer to as the \emph{relative space}.
This space is quasi-isometric to the \emph{relative} or
\emph{electrified Cayley graph}. The relative Cayley graph $\Ghat$ is
the graph formed by adding a vertex for each coset $gH_i$, and then
connecting this vertex to each element of the coset $gH_i$ by an edge
of length one-half.  However, we will find it convenient to work with
$\Grel$ rather than the electrified Cayley graph. We will write
$\nhat{x}$ for the relative length of $x$, and $\dhat{x,y}$ for the
relative distance from $x$ to $y$. A geodesic in the relative space is
a word of minimal length in the generators $A \cup \H$.

We say a finitely generated group $G$ is \emph{weakly relatively
  hyperbolic}, relative to a finite list of subgroups $\H$, if the
relative space $\Grel$ is $\delta$-hyperbolic.  This definition of
relative hyperbolicity is due to Farb \cite{farb}, and is more general
than the definition of strong relative hyperbolicity due to Gromov
\cite{gromov}, see also Bowditch \cite{bb2}. Osin \cite{osin} contains
a detailed discussion of several definitions of relative
hyperbolicity, and their relations.

We may consider the relative metric on the mapping class group with
respect to the following collection of subgroups.  Let $\{\alpha_1,
\ldots, \alpha_n\}$ be a list of representatives of orbits of simple
closed curves in $\S$, under the action of the mapping class group. We
may choose the $\a_i$ to be distance at most one from the curves
corresponding to the standard generators illustrated in Figure
\ref{picture21}. Let $H_i = \fix(\alpha_i$) be the subgroup of $G$
fixing $\alpha_i$. Masur and Minsky \cite{mm1} have shown that the
resulting relative space is quasi-isometric to the complex of
curves. As the complex of curves is $\delta$-hyperbolic, this shows
that the mapping class group is weakly relatively hyperbolic.

\begin{theorem}{\cite{mm1}*{Lemma 7.2}}
Let $\G$ be a Cayley graph for the mapping class group, and let
$\Ghat$ be the electrified Cayley graph with respect to subgroups
corresponding to stabilizers of representatives of orbits of the
mapping class group action on the vertices of the complex of
curves. Then the electrified Cayley graph $\Ghat$ is quasi-isometric
to the complex of curves $\C(\S)$.
\end{theorem}

In fact, the mapping class group is not strongly relatively
hyperbolic. There are a number of approaches to this available, of
which Karlsson and Noskov \cite{kn} seems to be the earliest, see also
Bowditch \cite{bb1}, Anderson, Aramayona and Shackleton \cite{aas},
and Behrstock, Drutu and Mosher \cite{bdm} for a useful discussion of
work in this area.

The quasi-isometry from the mapping class group with the relative
metric to the complex of curves can be defined by sending a group
element $g$ to $g(x_0)$, where $x_0$ is the basepoint of the complex
of curves. We will always measure distances in $\C(\S)$ between points
that lie in the zero-skeleton, so the complex of curves is effectively
a discrete space. Furthermore, the relative distance between two \mcg
elements is zero if and only if they are the same, so we only need an
additive quasi-isometry constant on the right hand side of the
inequalities below.
\begin{equation} \label{equation:qi} 
\frac{1}{\Kcc} \dhat{g,h} -\ccc \leqslant d_\C(gx_0,hx_0) \leqslant \Kcc \dhat{g,h}. 
\end{equation}

The \emph{Gromov boundary} of a $\delta$-hyperbolic space is the set
of endpoints of quasi-geodesic rays, where two quasi-geodesic rays are
said to be equivalent if they stay a bounded distance apart. The
Gromov boundary of a non-locally compact space need not be compact,
and in fact the Gromov boundary of the complex of curves is not
compact. As the complex of curves and the relative space $\Grel$ are
quasi-isometric, we may identify their Gromov boundaries. If $X$ is a
subset of either the complex of curves or the relative space, we will
write $\overline X$ for the union of $X$ and its limit points in the
Gromov boundary. The Gromov boundary of the complex of curves may be
described in terms of Teichm\"uller space and the Thurston
compactification of Teichm\"uller space.

\emph{Teichm\"uller space} can be defined as the space of complete
finite volume hyperbolic metrics on $\S$, modulo isometries isotopic
to the identity, which we shall denote $\T(\S)$. Topologically this is
homeomorphic to $\R^{6g+2p-6}$, where $g$ is the genus of $\S$, and
$p$ is the number of punctures.  Thurston showed that the space of
projective measured foliations $\PMF$ can be used to compactify
Teichm\"uller space in such a way that the action of the mapping class
group extends continuously to $\PMF$, which is a sphere of dimension
$6g+2p-7$.  The space of projective measured foliations $\PMF$ is
essentially the same as the space of projective measured laminations
$\PML$, as the two spaces are homeomorphic by a homeomorphism which
commutes with the action of the mapping class group.

Given $\e >0$, for each $\g \in \C_0(\S)$ define $\Thin(\g)$ to be the
subset of $\T(\S)$ in which the length of the geodesic representative
of $\g$ is at most $\e$. We may choose $\e$ to be sufficiently small
so that some collection $\Thin(\g_1), \ldots, \Thin(\g_n)$ has
non-empty intersection if and only if the curves $\g_1, \ldots, \g_n$
can be realised disjointly on the surface. Add a new point $x_\g$ for
each simple closed curve $\g$, and a new edge of length one-half
connecting $x_\g$ to each point in $\Thin(\g)$. The resulting space is
called \emph{electrified Teichm\"uller space}, $\T_{el}(\S)$. Masur
and Minsky \cite{mm1} show that $T_{el}(\S)$ is quasi-isometric to the
complex of curves $\C(\S)$, which is $\delta$-hyperbolic. In
particular, this means their Gromov boundaries are the same.

Klarreich \cite{klarreich}*{Theorem 1.1}, see also Hamenst\"adt
\cite{hamenstadt}, has identified the Gromov boundary of the complex
of curves with the space of minimal foliations on $\S$, with the
topology obtained from $\PMF$ by forgetting the measures. We shall
denote this space by $\fmin$, and the analogous space of laminations
$\lmin$.  A foliation is minimal if no leaf is a simple closed curve.
Minimal foliations correspond to laminations which contain no simple
closed curve, and which are not disjoint from any simple closed
curves, and we will call such laminations \emph{filling laminations}.
A \pA element has a unique pair of fixed points in $\PMF$ which are
uniquely ergodic foliations, and in particular they are minimal
foliations, and so a \pA element also has a unique pair of fixed
points in $\fmin$ or $\lmin$.

\subsection{Subgroups of the mapping class group}

We now summarize various well-known properties of subgroups of the
mapping class group that we will use in this paper. We provide
arguments for the sake of completeness, but these are either
elementary, or follow arguments from Ivanov \cite{ivanov}.

Given a subgroup $F$ of $G$ we will write $C(F)$ for the centralizer
of $F$, i.e. the subgroup of $G$ consisting of elements which commute
with elements of $F$, and we shall write $N(F)$ for the normalizer of
$F$, i.e. the subgroup of $G$ consisting of all elements $g$ such that
$gFg^{-1} = F$.  We will write $C(g)$ for the centralizer of the
cyclic subgroup generated by $g$. If $F$ is finite, then the
centralizer has finite index in the normalizer. A subgroup of the
mapping class group is \emph{non-elementary} if it contains a pair of
\pA elements with distinct pairs of fixed points in the Thurston
boundary $\PML$. A subgroup of the mapping class group is
\emph{reducible} if it preserves a finite collection of disjoint
simple closed curves.  Ivanov \cite{ivanov} showed that any infinite
subgroup of the mapping class group is either reducible or contains a
\pA element, and the centralizer of a \pA element $g$ is virtually
cyclic, and consists precisely of those elements which preserve the
fixed points of $g$.  Furthermore, if a subgroup consists entirely of
periodic elements, then the subgroup is finite, see also Birman,
Lubotzky and McCarthy, \cite{blm}, McCarthy \cite{mccarthy} and
McCarthy and Papadopoulos \cite{mp}.

It will be useful to know certain properties of finite subgroups of
the \mcg.  Let $F$ be a finite subgroup of the mapping class group. By
the Nielsen realization theorem, due to Kerckhoff \cite{kerckhoff},
there is a hyperbolic structure for $\S$ on which $F$ acts by
isometries. Elements of $F$ need not act freely on $\S$, so the
quotient $\O = \S /F$ may be an orbifold. Equivalently, in terms of
conformal structures, there is an $F$-invariant conformal structure on
$\S$ and the quotient map $\S \to \O$ is a branched cover.  This
covering is normal, and so gives rise to a map $\theta \colon \pi_1 \O
\to F$. The Teichm\"uller space of the quotient surface $\T(\O)$ is
the same as the Teichm\"uller space of the surface obtained by
replacing all branch points with punctures, and is isometrically
embedded in $\T(\S)$ as the fixed point set of $F$, i.e. $\{x \in
\T(\S) \mid f(x) = x \text{ for all } f \in F\}$, which we shall write
as $\fix_\T(F)$.  The mapping class group of the orbifold $\O$, which
we shall denote $G_\O$, consists of isotopy classes of homeomorphisms
from $\O$ to itself, which preserve the set of punctures, and also
preserve the sets of orbifold points or branch points of a given
index.

We say a homeomorphism of $\S$ is \emph{$F$-equivariant} if it
preserves pre-images of points under the covering map $\S \to \O$.  We
will write $G_F$ for the subgroup of $G$ consisting of elements with
$F$-equivariant representatives, which is also known as the
\emph{relative mapping class group} with respect to $F$. If an
$F$-equivariant homeomorphism of $\S$ is isotopic to the identity,
then it is isotopic by an $F$-equivariant isotopy, see Maclachlan and
Harvey \cite{mh}, Zieschang \cite{zieschang} and Birman and Hilden
\cite{bihi}, so this gives a well defined map $G_F \to G_\O$, with
kernel $F$. We now show that the image of $G_F$ in $G_\O$ has finite
index in $G_\O$.

\begin{proposition} \label{prop:subgroup image}
The image of the $F$-equivariant subgroup $G_F$ has finite index in
the mapping class group $G_\O$ of the quotient surface $\O$. 
\end{proposition}

\begin{proof}
A map $g \colon \O \to \O$ is covered by a
map $\widetilde g \colon \S \to \S$ if and only if $g_*(\ker \theta)
\subset \ker \theta$, where $\theta \colon \pi_1 \O \to \pi_1 \O$ is
the map induced by the normal covering. As $F$ is finite, $\ker
\theta$ is a finite index subgroup of $\pi_1 \O$, and there are only
finitely many finite index subgroups of that index, so the subgroup of
$G_\O$ which leaves $\ker \theta$ invariant has finite index in
$G_\O$. Therefore the image of $G_F$ has finite index in $G_\O$.
\end{proof}

We define the complex of curves of an orbifold $\C(\O)$ to be the
complex of curves of the surface obtained by treating the orbifold
points as punctures.  The curve complex $C(\O)$ is non-empty, unless
$\O$ is a triangle orbifold.  We will write $\fix_K(F)$ for the
\emph{coarse fixed set} of $F$ acting on the complex of curves
$\C(\S)$, i.e. all simple closed curves $x$ such that $\dc{x, fx} \le
K$, for all $f \in F$.  There is a map between the zero-skeletons
$\C_0(\O) \to \C_0(\S)$ which sends a simple closed curve in $\O$ to
its pre-image in $\C_0(\S)$, which is $F$-invariant and has diameter
at most one. Therefore this map is coarsely well defined, and sends
$\C_0(\O)$ to $\fix_1(F)$, which is empty if $\O$ is a triangle
orbifold.  We say a group acts \emph{coarsely transitively} on a
metric space $(X,d)$ if there is a constant $K$ such that the
$K$-neighbourhood of the orbit of a point in $X$ contains all of
$X$. We now show that the centralizer $C(F)$ acts coarsely
transitively on $\fix_1(F)$.

\begin{proposition} \label{prop:coarsely transitive}
Let $F$ be a finite subgroup of the \mcg. Then the centralizer $C(F)$
acts coarsely transitively on $\fix_1(F)$.
\end{proposition}

\begin{proof}
As $F$ is finite, the centralizer $C(F)$ has finite index in the
normalizer $N(F)$. By Proposition \ref{prop:subgroup image}, there is
a homomorphism with finite kernel from the normalizer $N(F)$ onto a
finite index subgroup of $G_{\O}$. The mapping class group $G_\O$ acts
coarsely transitively on $\C(\O)$, so $N(F)$ also acts coarsely
transitively on $\C(\O)$. The complex of curves $\C(\O)$ maps coarsely
onto $\fix_1(F)$, the set of all simple closed curves moved distance
at most one by all elements of $F$, and so $N(F)$, and hence $C(F)$,
acts coarsely transitively on $\fix_1(F)$.
\end{proof}

In the case that $\O$ is a triangle orbifold, the mapping class group
of $\O$, and hence the centralizer $C(F)$, is finite, and $\fix_1(F)$
is empty, so this statement is vacuously true.

The quasi-isometry from $\Grel$ to $\C(\S)$ is given by $g \mapsto
g(x_0)$, where $x_0$ is a fixed basepoint in $\C(\S)$.  As $C(F)$
preserves $\fix_1(F)$, the distance from the basepoint $x_0$ to
$\fix_1(F)$ is the same as the distance from $g(x_0)$ to $\fix_1(F)$
for any $g \in N(F)$. This implies that the image of $N(F)$ lies in a
bounded neighbourhood of $\fix_1(F)$, and by Proposition
\ref{prop:coarsely transitive}, the centralizer $C(F)$, and hence the
normalizer $N(F)$, act coarsely transitively on $\fix_1(F)$, so they
have the same limit sets in $\lmin$.  We now observe that the limit
set of $\fix_1(F)$, and hence the limit set of the centralizer $C(F)$,
consists of the filling laminations fixed by $F$. We will write
$\fix_\d(F)$ for the fixed set of $F$ in the Gromov boundary $\d
\C(\S)$.

\begin{proposition} \label{prop:centralizer limit}
If $F$ is a finite subgroup of the mapping class group of a
non-sporadic surface, then the limit set of the centralizer $C(F)$ is equal to
$\fix_\d(F)$, the fixed set of $F$ in $\lmin$.
\end{proposition}

\begin{proof}
If $\l$ is a filling lamination
in $\overline{\fix_1(F)}$, then it is $F$-invariant, as it is a limit of
$F$-invariant sets of disjoint simple closed curves. Conversely, if
$\l$ is a $F$-invariant lamination, then $\l / F$ is a lamination in
$\O$, so is the limit of a sequence of simple closed curves in $\O$.
The pre-images of these simple closed curves are $F$-invariant, and
hence lie in $\fix_1(F)$, so $\l$ is also a limit point of
$\fix_1(F)$. So $\overline{\fix_1(F)} \subset \d \C(\S)$ is equal to
the $F$-invariant filling laminations.
\end{proof}

The stabilizer of a fixed point set is the subgroup of $G$ which
leaves the fixed point set invariant, but does not necessarily fix it
pointwise. We now show that the stabilizer of these fixed sets is the
normalizer $N(F)$, which is also equal to the $F$-equivariant
subgroup.

\begin{proposition} \label{prop:equal subgroups}
Let $F$ be a finite subgroup of $G$. Then the following subgroups of
the mapping class group are equal: the $F$-equivariant subgroup $G_F$,
the normalizer $N(F)$, the stabilizer of $\fix_\T(F)$ and the
stabilizer of $\fix_\d(F)$.
\end{proposition}

\begin{proof}
Maclachlan and Harvey \cite{mh} and Birman and Hilden \cite{bihi}, show that the
$F$-equivariant subgroup of $G$ is equal to $N(F)$, the normalizer of
$F$, and furthermore, Maclachlan and Harvey \cite{mh} show that the
stabilizer of $\fix_T(F)$ in
$G$ is also equal to $N(F)$, as long as $\fix_\T(F)$ is non-empty,
which is always the case by the Nielsen realization theorem, due to
Kerckhoff \cite{kerckhoff}.  

The mapping class group $G$ acts continuously on both $\T(\S) \cup
\PML$ and $\C(\S) \cup \d \C(\S)$, and the Gromov boundary is the set
of filling laminations $\lmin$ which is a dense subset of $\PML$.  The
set of fixed points for $F$ in $\T(\S)$ is an isometrically embedded
copy of the Teichm\"uller space $\T(\O)$ for the quotient orbifold,
with limit set consisting of the laminations in $\PML$ fixed by
$F$. Therefore a mapping class group element preserves $\fix_\T(F)$ if
and only if it preserves $\fix_\d(F)$, so the stabilizer of
$\fix_\d(F)$ is equal to the stabilizer of $\fix_\T(F)$.
\end{proof}

In fact, the stabilizer of $\fix_\T(F)$ in the isometry group of
$\T(\S)$ is also equal to $N(F)$, as $\text{Isom}(\T(\S))$ is equal to
$G$ by Royden's theorem \cite{royden}, see Earle and Kra \cite{ek} for
the case of surfaces with punctures. 

\begin{proposition} \label{prop:small limit set}
Let $g$ be an element of the the \mcg of infinite order. Then the
limit set of the centralizer $C(g)$ in the Gromov boundary $\d \C(\S)$
consists of at most two points.
\end{proposition}

\begin{proof}
If $g$ is \pA, then the centralizer of $g$ is virtually
cyclic, and its image in the relative space $\Grel$ is a quasi-geodesic with endpoints
consisting of a pair of points in $\lmin$, namely the stable and
unstable laminations of $g$.  

If $g$ is reducible, then consider a power of $g$ which is pure, i.e.
there is a collection of disjoint simple closed curves $a_i$, each of
which is fixed by $g^n$, and $g^n$ acts on each complementary
subsurface as either the identity or as a \pA. As the centralizer
$C(g)$ is contained in the centralizer $C(g^n)$, it suffices to
consider the case in which $g$ is pure. We now show that the fixed
point set of a reducible pure element $g$ has bounded diameter in the
relative metric. If there is a subsurface on which $g$ acts as a \pA,
then each curve distance two from any $a_i$ has infinitely many images
under powers of $g$, so the fixed set of $g$ in $\C(\S)$ has bounded
diameter. If $g$ acts as the identity on each complementary
subsurface, then $g$ must act as a power of a Dehn twist on one of the
fixed curves $a_i$. Again, any curve distance two or more from $a_i$
has infinitely many images under powers of $g$, and so the fixed set
of $g$ has bounded diameter in this case as well.  The fixed set of
$g$ is preserved by $C(g)$, as if $x \in \fix_0(g)$, and $h \in C(g)$
then $gh(x) = hg(x) = h(x)$, so $h(x) \in \fix_0(g)$.  Furthermore,
the image of $C(g)$ in the relative space $\Grel$ also has bounded
diameter, as the image of $h$ in $\C(\S)$ is $h(x_0)$. As the fixed
set $\fix_0(g)$ is preserved by $h$, the distance from $h(x_0)$ to the
fixed set is the same as the distance from $x_0$ to the fixed set, so
$\dc{h(x_0), \fix_0(g)}$ is independent of $h \in C(g)$. This
implies that the diameter of $C(g)$ is bounded in the relative metric,
and so the limit set of $C(g)$ in the Gromov boundary $\lmin$ is
empty.
\end{proof}

This implies that if $H$ is a non-elementary subgroup of the mapping
class group, and its limit set is contained in the centralizer $C(F)$
of some subgroup $F$ of the \mcg, then every element of $F$ is finite,
and so $F$ is a finite subgroup of the \mcg.

\begin{proposition} \label{prop:limit set}
Let $H$ be a non-elementary subgroup whose limit set is contained in the
limit set of the centralizer of a finite subgroup $F$ in
the \mcg $G$. Then $H$ is contained in the normalizer $N(F)$. 
\end{proposition}

\begin{proof}
We may assume that $F$ is a maximal subgroup such that $\overline H
\subset \overline{C(F)}$, as $F$ is finite, and there is an upper
bound on the size of any finite subgroup of $G$, which depends on the
surface $\S$. The \mcg acts continuously on $\C(\S) \cup \d \C(\S)$,
and $\overline{H}$ is $H$-invariant, therefore $h \overline{H} \subset
h \overline{C(F)}$ for any $h \in H$. In particular,
\[ \overline{H} \subset \bigcap_{h \in H} h\overline{C(F)}. \] 
A translate $h \overline{C(F)}$ is equal to $\overline{C(hFh^{-1})}$,
and the intersection of two centralizers $C(F_1)$ and $C(F_2)$ is
equal to the centralizer of the group generated by $F_1 \cup
F_2$. Therefore $\overline{H}$ is contained in $\overline{C(F')}$,
where $F'$ is the group generated by all conjugates of $F$ by elements
of $H$. The limit set of $\overline{C(F')}$ is the intersection of the
limit sets $\overline{C(f)}$ of each element $f \in F'$. If any
element $f$ in $F'$ were infinite order, then by Proposition
\ref{prop:small limit set}, the limit set $\overline{C(F)}$ would
consist of at most two points. As $H$ is non-elementary, its limit set
$\overline{H}$ contains infinitely many points, so all elements of
$F'$ are finite order, and so $F'$ is a finite subgroup of the \mcg
$G$. However, $F'$ contains $F$, so by the maximality of $F$, the
subgroup $F'$ is in fact equal to $F$. This implies that $F$ is
invariant under conjugation by elements of $H$, i.e. $H$ is contained
in the normalizer $N(F)$, as required.
\end{proof}

Let $H^+$ be a semi-group in $G$, which generates a non-elementary
subgroup of the mapping class group, and let $F$ be a finite subgroup
of the \mcg.  We now show that either $\overline{C(F)}$ has infinitely
many distinct, though not necessarily disjoint, images in $\lmin$
under $H^+$, or else $H^+$ is contained in the normalizer of $F$.

\begin{proposition} \label{prop:infinite images}
Let $H^+$ be a semi-group in the \mcg which generates a non-elementary
subgroup $H$. Then for any finite subgroup $F$ of the \mcg then either
there are infinitely many images of $\overline{C(F)}$ under $H^+$, or
else $H$ is contained in the normalizer $N(F)$.
\end{proposition}

\begin{proof}
Let $H^+$ be a semi-group in $G$ which generates a non-elementary
subgroup. Suppose there is a finite group $F$ such that there are only
finitely many images of $\overline{C(F)}$ under $H^+$. There are only
finitely many conjugacy classes of finite subgroups in the mapping
class group, so we may assume we have chosen a maximal finite subgroup
$F$ with this property.  Each element of $H^+$ acts as a finite
permutation on the $H^+$-orbit of $\overline{C(F)}$, which implies
that the inverse of each element in $H^+$ also acts as a finite
permutation on the $H^+$-orbit of $\overline{C(F)}$. As $H$ is
generated by $H^+$, together with the inverses of elements in $H^+$,
the entire group $H$ acts as finite permutations on the $H^+$-orbit of
$\overline{C(F)}$.  This implies that $H$ contains a finite index
subgroup $H'$ that preserves $\overline{C(F)}$. As the stabilizer of
$\overline{C(F)}$ is the normalizer $N(F)$, by Proposition
\ref{prop:equal subgroups}, this implies that $H'$ is a subgroup of
$N(F)$.  Therefore the limit set of $H'$ is contained in the limit set
of $N(F)$. As the limit set of a finite index subgroup is the same as
the limit set of the original group, this implies that the limit set
of $H$ is contained in the limit set of $C(F)$.  Therefore, by
Proposition \ref{prop:limit set}, $H$ is contained in the normalizer
$N(F)$.
\end{proof}

\section{Relative conjugacy bounds}
\label{section:conjugacy}

In this section we show that if two elements of the mapping class
group $G$ are conjugate, then the relative length of the shortest
conjugating word is bounded in terms of the relative lengths of the
two conjugate elements.

\begin{theorem} 
\label{theorem:conjugacy} 
Let $a$ and $b$ be conjugate elements of the mapping class group of a
non-sporadic surface.  Then there is a conjugating word $w$ of
relative length $\nhat{w} \leqslant \Kconj(\nhat{a}+\nhat{b})$, for
some constant $\Kconj$ which only depends on the surface $\S$.
\end{theorem}

The proof of this result relies on the fact that there is a
``discreteness'' for the action of the mapping class group on the
complex of curves. This means that pseudo-Anosov elements, which act
on the complex of curves as hyperbolic isometries, have a minimal
translation length, and reducible or periodic isometries, which act as
elliptic isometries, act coarsely transitively on their coarse fixed
sets.

We will prove this by considering each of the three different types of
elements of the mapping class group in turn. Each case will produce a
different constant, but we will then choose $\Kconj$ to be the maximum
of the three constants. Without loss of generality, we may choose $K$
to be at least $1$, and we shall do this, as it will enable us to
simplify expressions in subsequent sections.

\subsection{Pseudo-Anosov elements}

Masur and Minsky have shown that the mapping class group has
(non-relative) conjugacy bounds for pseudo-Anosov elements, i.e.  if a
$a$ and $b$ are conjugate pseudo-Anosov elements, then there is a
conjugating element $w$ such that $\norm{w} \leqslant
K(\norm{a}+\norm{b})$. Although, a priori, the properties of having
relative or non-relative conjugacy bounds for a weakly relatively
hyperbolic group may be independent, our argument for the
pseudo-Anosov case is clearly modelled on Masur and Minsky's argument
from \cite{mm2}, and in fact the relative case is substantially
simpler.

Pseudo-Anosov elements act on the complex of curves as hyperbolic
isometries. A hyperbolic isometry $h$ of a $\delta$-hyperbolic space
has a \emph{quasi-axis}, which is a bi-infinite quasi-geodesic
$\axis$, such that $\axis$ and $h^k\axis$ are $2\delta$ fellow
travellers for all $k$. In the case of the complex of curves, Masur
and Minsky \cite{mm2} have shown that we may choose the quasi-axis to
be a geodesic, so we will do this, and we will refer to it as an
\emph{axis} for $h$.

The following theorem of Masur and Minsky \cite{mm1} shows that there
is a lower bound on the translation length of a pseudo-Anosov element
acting on the complex of curves, which only depends on the surface $\S$.

\begin{theorem}{\cite{mm1}*{Proposition 3.6}}
\label{theorem:translation distance}
Let $h$ be a pseudo-Anosov element of the mapping class group of a
surface which is not a sphere with three or fewer punctures. Then
there is a constant $c$ such that $\dc{x,h^nx} \geqslant c|n|$, for
all $x \in \C(\S)$. The constant $c$ depends on the surface $\S$, but
is independent of the pseudo-Anosov element $h$.
\end{theorem}

In \cite{mm1}, Masur and Minsky state this for non-sporadic surfaces,
but the result also holds in the case of a torus with one or fewer
punctures or a four-punctured sphere.

We first show there is a lower bound for the distance a point $x$ in
$\C(\S)$ is moved by $h$, in terms of the translation distance of $h$
along its axis $\a$, and the distance of $x$ from the axis.

\begin{lemma} \label{lemma:axis}
Let $\axis$ be an axis for $h$, and let $x$ be an element of the complex
of curves, then $\dc{x,\axis} \le \K1 \dc{x,hx}$, where
$\K1 = 4\delta/c$.
\end{lemma}

\begin{proof}
Let $y$ be the closest point to $x$ on the axis $\axis$. Choose an
integer $n \geqslant 8\delta/c$, which only depends on the surface $\S$, such that
$d(y,h^n y) \geqslant 8\delta$, where $\delta$ is the
$\delta$-hyperbolicity constant for $\C(\S)$. Consider a geodesic
$[x,h^n x]$, which forms a side of a quadrilateral, together with
$[x,y], [y,h^n y]$ and $[h^n x,h^n y]$.  This is illustrated in Figure
\ref{picture33} below.

\begin{figure}[H]
\begin{center}
\epsfig{file=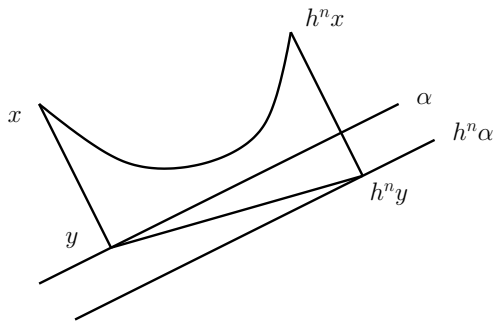, height=120pt}
\end{center} 
\caption{Translation along an axis.}\label{picture33}
\end{figure}

In a $\delta$-hyperbolic space any side of a quadrilateral is
contained in a $2\delta$-neighbourhood of the other three. As $[x,y]$
and $[h^nx,h^ny]$ are shortest paths to $\axis$ and $h^n\axis$, which
lie in $2\delta$ neighbourhoods of each other, any point on $[y,h^ny]$
which is at least $4\delta$ away from each of its endpoints must also
be at least $2\delta$ away from $[x,y]$ and $[h^nx,h^ny]$, so must be
$2\delta$ close to $[x,h^nx]$.  So the length of $[x,h^nx]$ must be at
least $2\dc{x,y} + \dc{y,h^ny} - 8\delta \geqslant 2\dc{x,y}$.

As $\dc{x,h^nx} \le n\dc{x,hx}$, this means $\dc{x,y} \le
\frac{n}{2}\dc{x,hx}$, so the claim follows, with $\K1 = n/2 =
4\delta/c$, which only depends on the surface $\S$, as required.
\end{proof}

Now let $a$ and $b$ be two conjugate elements of $G$, and let $w$ be
some conjugating element, such that $a=wbw^{-1}$. Let $\a$ be an axis
for $a$ and let $\b$ be an axis for $b$. Let $y$ be the closest
point in $\a$ to $x_0$, and let $z$ be the closest point on
$\b$ to $x_0$, as illustrated below in Figure \ref{picture28}.

\begin{figure}[H]
\begin{center}
\epsfig{file=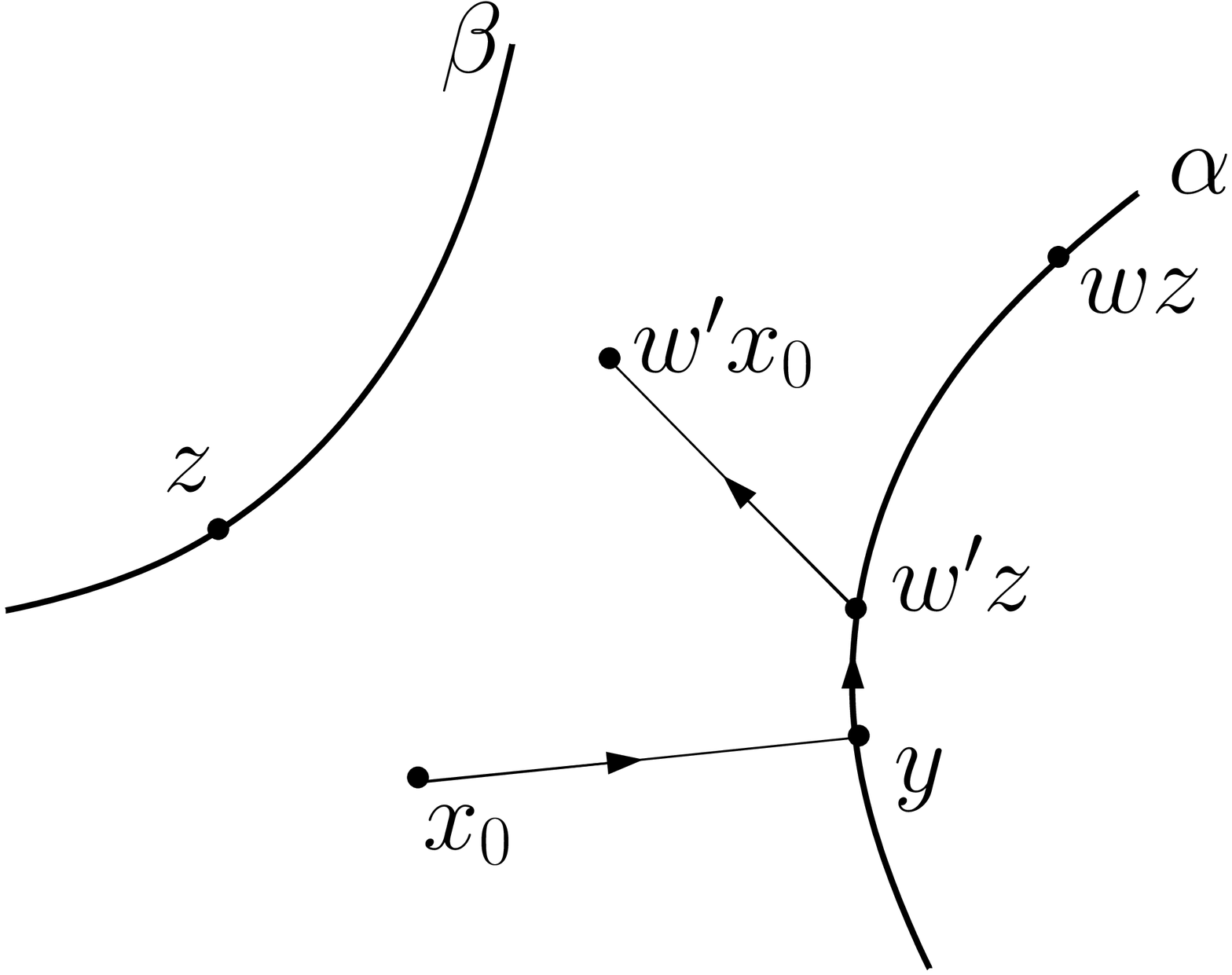, height=120pt}
\end{center} 
\caption{Estimating $\nhat{w'}$.}\label{picture28}
\end{figure}

The spaces $\Grel$ and $\C(\S)$ are quasi-isometric, so
$\frac{1}{\Kcc} \nhat{a} \leqslant \dc{x_0,ax_0} \leqslant \Kcc
\nhat{a} $, where $\Kcc$ is the quasi-isometry constant from
\eqref{equation:qi}.  Therefore, using Lemma \ref{lemma:axis} above,
$\dc{x_0,y} \leqslant \Kcc \K1 \nhat{a}$, and $\dc{x_0,z} \leqslant
\Kcc \K1 \nhat{b}$. The conjugating element $w$ takes the axis
$\axisb$ to the axis $\axisa$, so $wz \in \axisa$. The translation
length of $a$ is at most $\dc{x_0,ax_0} \leqslant \Kcc \nhat{a}$, and
so $a$ acts $\Kcc \nhat{a}$-coarsely transitively on its axis $\a$, so
there is a $k$ such that $\dc{y,a^kw z } \leqslant \Kcc \nhat{a}$. We
can choose to conjugate $b$ to $a$ by $w'=a^kw$ instead of $w$, and we
now show that the relative length of $w'$ is bounded.

The relative length of $w'$ is bounded in terms of the distance it
moves the base point,
\begin{align*}
\frac{1}{\Kcc} \nhat{w'} - \ccc & \leqslant \dc{x_0, w' x_0}. \\
\intertext{But this distance is at most the length of the path formed
  by going from $x_0$ to $y \in \axisa$, then along the axis to $w'
  z$, and then to $w' x_0$, as illustrated above in Figure
  \ref{picture28}.}
\frac{1}{\Kcc} \nhat{w'} - \ccc & \leqslant   \dc{x_0,y} +  \dc{y,w' z } +  \dc{w'z ,w' x_0 } \\
\intertext{The first term is roughly $\nhat{a}$, the second term is at
  most the translation length of $a$, and the final term is equal to
  $\dc{z,x_0}$, which is roughly $\nhat{b}$.}  
\frac{1}{\Kcc} \nhat{w'} - \ccc & \leqslant \Kcc \K1 \nhat{a} + \Kcc \nhat{a} + \Kcc \K1 \nhat{b}
\end{align*}
The only mapping class group element of relative length zero is the
identity, so we may assume that both $\nhat{a}$ and $\nhat{b}$ are at
least $1$. Therefore $\nhat{w'} \leqslant \Kcc^2(\K1 + 1)\nhat{a} +
(\Kcc^2 \K1 + q )\nhat{b} \leqslant \Kcc^2(\K1 + 1 + q)(\nhat{a} +
\nhat{b})$.  So we have shown that for pseudo-Anosov elements we may
choose the conjugacy bound constant to be $\Kconj = 2\Kcc^2(\K1 + 1 +
q)$, which only depends on the surface $\S$.

\subsection{Reducible elements}

Reducible elements of the mapping class group act on the complex of
curves as elliptic isometries. A reducible element $h$ leaves
invariant a collection of disjoint simple closed curves. Following
Ivanov \cite{ivanov}, we say an element of the mapping class group is
\emph{pure}, if there is a disjoint collection of simple closed curves
$\s(h)$ which are fixed individually by $h$, such that each
complementary component of $\s(h)$ is also fixed, and furthermore $h$
acts on each complementary component as either a pseudo-Anosov element
or the identity.  If the collection of simple closed curves $\s(h)$
has the property that no simple closed curve with non-zero
intersection number with $\s(h)$ is fixed by $h$, then $\s(h)$ is
called a \emph{canonical reduction system} for $h$.  If $h$ is not
pure, then we define the canonical reduction set $\s(h)$ to be the
canonical reduction set of some pure power of $h$.  Given a reducible
element $h$ of the mapping class group, we can raise $h$ to some power
$m_1$ so that $h^{m_1}$ is pure.  There is a power $m_1$ which works
for all reducible elements of the mapping class group, depending only
on $\S$.

We start by showing that there is a lower bound on the distance a
point is moved by $h$ in terms of its distance from the fixed curves.

\begin{lemma} 
\label{lemma:reducible}
Let $h$ be a reducible element of the mapping class group of a
non-sporadic surface, with canonical reduction set $\s(h)$, and let
$y$ be a vertex of the complex of curves. Then $\dc{\s(h),y} \leqslant
\K2\dc{y,hy} + 3$, for some constant $\K2$, which only depends on the
surface $\S$.
\end{lemma}

\begin{proof}
We may raise $h$ to some power $m_1$, depending on $\S$, such that
$h^{m_1}$ is pure, and we shall just write $h$ for $h^{m_1}$ from now on.

Given a connected subsurface $\S' \subset \S$ with essential boundary
components, which is not a three-punctured sphere, Masur and Minsky
\cite{mm2} define a \emph{subsurface projection} $\pi:\C(\S) \to
\C(\S) \cup \varnothing$, which we now describe. Given a simple closed
curve $y$ in $\C(\S)$ which intersects $\d \S'$ efficiently, we can
take a component of $y \cap \S'$ and complete it to a simple closed
curve in $\S'$ by adding a boundary parallel arc. This actually gives
a map from $\C(\S)$ to subsets of $\C(\S')$, but the image has bounded
diameter independent of $y$, so this is coarsely well-defined. If $y$
does not intersect $\S'$ then we send it to the empty set. In the case
where the subsurface is an annulus, Masur and Minsky provide an
appropriate definition of a complex for the annulus, which is coarsely
equivalent to $\Z$, and which roughly speaking counts how many times
the arcs of $y$ wrap around the annulus. We will not give the details
of this definition here, as the only property we will use is that a
Dehn twist in the core curve of the annulus acts with translation
distance one on the complex for the annulus.

We will use the following bounded geodesic image theorem of Masur and
Minsky \cite{mm2}, which says that a geodesic which is not close to
$\d \S' \subset \C(\S)$ projects to a bounded set under the subsurface
projection.

\begin{theorem} \cite{mm2}*{Theorem 3.1} \label{theorem:bounded image}
Let $\S'$ be an essential connected subsurface of $\S$, which is not
a three-punctured sphere, and let $\g$ be a geodesic segment in
$\C(\S)$, such that $\pi(v) \not = \varnothing$ for every vertex $v
\in \g$. Then there is a constant $M$, which only depends on $\S$,
such that the diameter of $\pi(\g)$ is at most $M$.
\end{theorem}

Suppose there is a component $\S'$ of $\S \setminus \s(h)$ on which
$h$ acts as a pseudo-Anosov element. Then by Theorem
\ref{theorem:translation distance}, the translation distance of the
pseudo-Anosov element $h|_{\S'}$ on the complex of curves $\C(\S')$ is
bounded below, so we can raise $h$ to some power $m_2$, which only
depends on $\S$, so that $h^{m_2}$ has translation distance at least
$M+1$ on $\C(\S')$. Otherwise, if there are no components where $h$
acts as a \pA, then there is a simple closed curve $x$ in $\s(h)$ such
that $h$ acts as a Dehn twist in a neighbourhood of $x$, then choose
$\S'$ to be a regular neighbourhood of $x$, and then $h^{M+1}$ has
translation distance $M+1$ on $\C(\S')$.

We may assume $h$ acts as the identity in a neighbourhood of the
boundary of the surface $\S$. Let $y$ be a simple closed curve in
$\C(\S)$, then the subsurface projection map $\pi:\C(\S) \to \C(\S')$
only alters $y \cap \S'$ in a neighbourhood of the boundary of the
subsurface $\S'$, so the subsurface projection $\pi$ and $h|_{\S'}$ commute.

Now let $y$ be a vertex of the complex of curves such that $\dc{\s(h),y}
\geqslant 3$, and let $\g$ be a geodesic in $\C(\S)$ from $y$ to
$hy$. If $\dc{\g,\s(h)} \geqslant 3$, then Theorem \ref{theorem:bounded
  image} implies that the image of the projection $\pi(\g)$ in
$\C(\S)$ has diameter at most $M$. However, $\pi(\g)$ contains both
$\pi(y)$ and $\pi(hy) = h(\pi(y))$, which are distance at least
$M+1$ apart, which gives a contradiction. So $\g$ must pass within
distance three of $\s(h)$, which implies that $\dc{y, hy}$ is at least
$2\dc{\s(h),y} - 6$.

As we may have raised $h$ to some power, we may take $\K2 =
\frac{1}{2} m_1 \max\{ m_2, M+1 \} +6$, which only depends on
$\S$. Furthermore, we assumed $\dc{\g, \s(h)} \ge 3$, so we also need
an additive term of $3$ in the inequality.
\end{proof}

Now suppose that $a$ and $b$ are conjugate reducible elements of the
mapping class group, with $a = wbw^{-1}$, for some conjugating word
$w$. Let $A$ be the canonical reduction set for $a$, and let $B$ be
the canonical reduction set for $b$.

By Lemma \ref{lemma:reducible} above, $\dc{x_0,A} \leqslant \K2
\dc{x_0, ax_0} + 3$, which in turn is at most $\Kcc \K2 \nhat{a} + 3$,
as $\Grel$ and $\C(\S)$ are $\Kcc$-quasi-isometric. Similarly,
$\dc{x_0,B} \leqslant \Kcc \K2 \nhat{b}$.  The conjugating element $w$
takes $B$ to $A$, and the diameter of the fixed sets is $2$, so
$\dc{x_0,wx_0}$ is at most the length of a path from $x_0$ to $A$, and
then from $A$ to $wx_0$, plus $2$. This implies
\begin{align*}
\dc{x_0,wx_0} & \leqslant \dc{x_0, A} + \dc{A, w x_0} + 2. \\ 
\intertext{The first term is bounded in terms of $\nhat{a}$, and the
  second term is equal to $\dc{B,x_0}$, as $wB = A$, and so is bounded
in terms of $\nhat{b}$. Therefore we obtain a bound in terms of the relative
lengths of $a$ and $b$,}
\dc{x_0, wx_0} & \leqslant \Kcc \K2 \nhat{a} + \Kcc \K2 \nhat{b} + 5. \\
\end{align*}
This means that $\nhat{w} \leqslant \Kcc^2 \K2 (\nhat{a}+\nhat{b}) +
\Kcc(5 + \ccc)$. We may assume $\nhat{a}$ and $\nhat{b}$ have length
at least one, so we may choose $K = \Kcc^2 \K2 + \Kcc(5 + \ccc)$,
which only depends on the surface $\S$.

\subsection{Periodic elements} \label{section:periodic}

Periodic elements of the mapping class group act on the complex of
curves as elliptic isometries. We now review some useful properties of
elliptic isometries, as described for example in Bridson and Haefliger
\cite{bh}.  In coarse geometry the analogue of the fixed point set of
an elliptic isometry of hyperbolic space is the \emph{$\e$-fixed set}
of a periodic element $h$, which is all points moved at most $\e$ by
the isometry, i.e. $\fix_{\e}(h) = \{x \in \C(\S) \mid \dc{x,hx}
\leqslant \e \}$.  Let $X$ be a bounded set in a $\delta$-hyperbolic
space $Y$, with radius $\rho = \inf \{ \rho \mid X \subset B_\rho(x)
\text{ for some } x \in X \}$, where $B_\rho(x)$ is the ball of radius
$\rho$ with center $x$. For $\e > 0$ the \emph{quasi-center} of $X$ is
$\qcenter_\e(X) = \{ y \in Y \mid X \subset B_{\rho+\e}(y)\}$.

\begin{lemma}{\cite{bh}*{Lemma 3.3}}
Let $X$ be a bounded set with quasi-center $\qcenter_\e(X)$. Then the diameter
of $\qcenter_\e(X)$ is at most $4\delta + 2\e$.
\end{lemma}

We shall choose $\e = \delta$, and define the \emph{coarse fixed
  set} of a periodic element $h$ to be $\fix_{6\delta}(h) = \{x \in
\C(\S) \mid \dc{x,hx} \leqslant 6\delta\}$. The orbit of any point
under $h$ is finite, and hence has a quasi-center, which is contained
in the coarse fixed set. In particular, the coarse fixed set is
non-empty.

\begin{lemma} 
\label{lemma:coarse}
Let $H$ be a finite cyclic subgroup of the \mcg of a non-sporadic
surface, generated by a periodic element $h$. Then the centralizer $C(H)$ acts $\K3$-coarsely
transitively on the coarse fixed set $\fix_{6\delta}(H)$. Furthermore,
for any $x$, $\dc{x,\fix_{6\delta}(h)} \leqslant \K3 \dc{x,hx}$. The
constant $\K3$ only depends on the surface $\S$.
\end{lemma}

\begin{proof}
Let $h$ be a periodic element of the mapping class group, with period
$n$, and let $H$ be the finite cyclic subgroup generated by $h$.  The
centralizer $C(H)$ acts coarsely transitively on the $1$-fixed set
$\fix_1(H)$, by Proposition \ref{prop:coarsely transitive}. We will
show that $\fix_1(H)$ and $\fix_K(H)$ have the same limit sets and are
quasiconvex, which implies that $C(H)$ acts coarsely transitively on
$\fix_K(H)$. The fact that the constant in the final inequality only depends on $\S$
then follows from the fact that there is an upper bound on the order
of a periodic element of the \mcg $G$, depending only the surface
$\S$.

We now show that $\fix_1(H)$ is a quasi-convex subset of $\C(\S)$. The
Teichm\"uller space $\T(\O)$, where $\O$ is $\S/ H$, can be identified
with the fixed set of $H$ in $\T(\S)$, which is isometrically
embedded, so Teichm\"uller geodesics between points in the fixed set
of $h$ in $\T(\S)$ are in fact contained in the fixed set of $H$ in
$\T(\S)$. Masur and Minsky \cite{mm1} show that Teichm\"uller
geodesics give rise to unparameterized quasi-geodesics in the complex
of curves, where an unparameterized quasi-geodesic is contained in a
uniform neighbourhood of a geodesic, and makes coarsely monotone, but
not necessarily coarsely uniform, progress along the geodesic.  This
means that the map from $\C(\O)$ to $\C(\S)$ is a quasiconvex
embedding, so the convex hull of $\overline{\fix_1(H)}$ is contained
in a bounded neighbourhood of $\fix_1(H)$.  In fact, Rafi and
Schleimer \cite{rs} have shown that this map is a quasi-isometric
embedding.

We now show that $\dc{x,\fix_{6\delta}(h)} \leqslant \K3 \dc{x,hx}$,
for any $x \in C(\S)$. Let $H . x$ be the orbit of $x$ under $H$, and
let $\qcenter_{\delta}(H . x)$ be the quasi-center of $H . x$. Suppose
$y \in \qcenter_{\delta}(H . x)$, then $H.x \subset
B_{\rho+\delta}(y)$, where $\rho$ is the radius of $H.x$, so in
particular, 
\[ \dc{x,\qcenter_{\delta}(H . x)} \le \rho + \delta. \]
The set $\qcenter_{\delta}(H . x)$ is $H$-invariant, as $H.x$ is
$H$-invariant, and has diameter at most $6\delta$, so
$\qcenter_{\delta}(H . x) \subset \fix_{6\delta}(H)$.  The diameter of
$H . x$ is at most $n \dc{x, h x}$, and so $ \dc{x,
  \qcenter_{\delta}(H . x)} \le n \dc{x, h x} + \delta$.  As
$\qcenter_\delta(H . x) \subset \fix_{6 \delta}(H)$, this implies
\[ \dc{x,\fix_{6\delta}(H)} \leqslant N \dc{x,hx} + \delta, \]
where $N$ is the maximum period of any periodic element of $G$, which
only depends on the surface $\S$. If $\dc{x, hx}$ is zero, then $x$ is
fixed by $h$, so $\dc{x, \fix_{6 \delta} (H)}$ is also zero, so we may
choose the constant $K_3$ here to be $N + \delta$.

Finally, we show that $\fix_{K}(H)$ is quasiconvex for $K \geqslant
6\delta$, and $C(H)$ acts coarsely transitively on $\fix_K(H)$. Let
$a$ and $b$ be points in $\fix_{K}(H)$, then the orbits $H . a$ and $H
. b$ each have diameter at most $n K$, and so the geodesics $[h^i a,
h^i b]$ are $2 \delta$-fellow travellers, outside of
$nK$-neighbourhoods of their endpoints, where $n$ is the period of
$h$.  Therefore, for $x \in [a, b]$, the orbit $H . x$ has diameter at
most $n(K + 2 \delta)$, and so $\dc{x, \qcenter_{\delta}(H . x)} \le n
(K + 2 \delta)+ \delta$. As $\qcenter_{\delta}(H . x) \subset
\fix_{6\delta}(H)$, which in turn is contained in $\fix_{K}(H)$, this
implies that every geodesic $[a, b]$ with endpoints in $\fix_K(H)$ is
contained in a $(N K + 2N \delta + \delta)$-neighbourhood of
$\fix_{K}(H)$, where $N$ is the largest order of a periodic element in
$G$, which depends only on $\S$. As $\fix_K(H)$ is quasiconvex, and
has the same limit set as $\fix_1(H)$, $\fix_K(H)$ lies in a bounded
neighbourhood of $\fix_1(H)$, and so as $C(H)$ acts coarsely
transitively on $\fix_1(H)$, $C(H)$ also acts coarsely transitively on
$\fix_K(H)$.  Therefore, we make take $K_3$ to be the maximum of $N +
\delta$, and the constant by which $C(H)$ acts coarsely transitively
on $\fix_{6\delta}(H)$.
\end{proof}

Let $a$ and $b$ be conjugate periodic elements of the mapping class
group, so $a = wbw^{-1}$, for some conjugating word $w$. Let $A$ be
the cyclic subgroup generated by $a$, and let $B$ be the cyclic
subgroup generated by $b$.  Let $y$ be the closest point in $\cfix{a}$
to $x_0$, and let $z$ be the closest point in $\cfix{b}$ to $x_0$.
The conjugating element $w$ takes $\cfix{b}$ to $\cfix{a}$, so $wz \in
\cfix{a}$. As $C(a)$ acts $\K3$-coarsely transitively on $\cfix{a}$,
there is $c \in C(a)$ such that $\dc{y,cw z} \leqslant \K3$. As $c$ is
in the centralizer of $a$, the element $cw$, which we shall denote
$w'$, conjugates $b$ to $a$.

We can estimate the relative length of the conjugating element
$\nhat{w'}$ in terms of the distance in the complex of curves from
$x_0$ to $w' x_0$, which is at most the length of the path from $x_0$
to $y$, then from $y$ to $w' z $, and finally from $w' z$ to $w'
x_0$. Figure \ref{picture28} also illustrates this case, if the axes
$\a$ and $\b$ are replaced with the coarse fixed sets $\cfix{a}$ and
$\cfix{b}$ respectively.
\begin{align*}
\dc{x_0,w' x_0} & \leqslant \dc{x_0, y} + \dc{y, w' z } + \dc{w' z
  ,w' x_0 } \\
\intertext{By Lemma \ref{lemma:coarse}, the first term on the right
  hand side is at most $\K3 \dc{x_0,ax_0}$. Similarly the final term
  on the right hand side is at most $\K3 \dc{x_0,bx_0}$. The middle
  term is bounded by $\K3$.}
\dc{x_0,w' x_0} & \leqslant \K3 \dc{x_0,a x_0} + \K3 + \K3
\dc{x_0,b x_0} \\
\intertext{The first term on the right hand side is at most $\Kcc \K3
  \nhat{a}$, using the quasi-isometry between $\Grel$ and
  $\C(\S)$. Similarly, the final term is at most $\Kcc \K3 \nhat{b}$.}
\dc{x_0,w' x_0} & \leqslant \Kcc \K3 \nhat{a} + \K3 + \Kcc \K3 \nhat{b}
\end{align*}
Therefore, $\nhat{w'} \leqslant \Kcc^2 \K3(\nhat{a}+\nhat{b}) + \Kcc
(\K3 + \ccc)$, and as we may assume that $a$ and $b$ have relative
length at least $1$, we may take $K = \Kcc^2 \K3 + \Kcc (\K3 + \ccc
)$, which only depends on the surface $\S$.

This completes the proof of Theorem \ref{theorem:conjugacy}.  The
conjugacy bound constants we have obtained in each of the above
sections may be different, however, we may choose $K$ to be the
maximum such constant for the three types of elements of the mapping
class group.

\section{Conjugates of relatively short elements} \label{section:coarse}

In this section we consider collections of elements which are
conjugates of relatively short elements, and we investigate how they
are distributed inside $G$. The results of this section hold for
any group $G$ which is weakly relatively hyperbolic and which has
relative conjugacy bounds. We shall refer to the constant of
hyperbolicity $\delta$ and the relative conjugacy bound constant $K$
as the \emph{group constants}. In particular, the results of this
section apply to the mapping class group, and the collection of
non-\pA elements in the \mcg consists of elements which are all
conjugate to elements of bounded relative length.

Let $R$ be a collection of elements of $G$, which are conjugate to
elements of relative length at most $B$. In some parts of $G$,
elements of $R$ are close together, in other parts of $G$ they are far
apart. We quantify this by defining $R_k$ to be the \emph{$k$-dense}
subset of $R$, consisting of all elements of $R$ which are distance at
most $k$ in $G$ from some other element of $R$, i.e.  $R_k=\{ r \in R
\mid \text{there is an } r' \in R \text{ with } r \not = r' \text{ and
} d(r,r') \leqslant k \}$. This definition uses word length in $G$, not
relative length.

Any element in $R_k$ differs from another element of $R$ by an element
$g$ of word length at most $k$, so in fact $R_k$ is the finite union
of sets $R \cap Rg$, as $g$ runs over all group elements of word
length at most $k$. The limit set of $R$, and hence of $Rg$, may be
the entire boundary, and this is the case, for example, if $R$ consists
of all non-\pA elements. However, the limit set of the intersection
$\overline{R \cap Rg}$ may be smaller than than the intersection of
the limit sets $\overline{R} \cap \overline{Rg}$. For example, if $R$
and $Rg$ are disjoint, then the limit set of their intersection will
be empty. In this section we will show that the elements of $R \cap
Rg$ are contained in a particular neighbourhood of the centralizer of
$g$, which we shall call a \emph{horoball neighbourhood}, as its
definition is reminiscent of the definition of a horoball in
hyperbolic space.  Let $X$ be a subset of $\Grel$, and let $L$ be a
constant. We define an $L$-horoball neighbourhood of $X$, which we
shall denote $\horo{X}{L}$, to be the union of balls in $\Grel$
centered at $x \in X$, of radius $\nhat{x}+L$, i.e.
\[ \horo{X}{L} = \bigcup_{x \in X} \widehat B_{\nhat{x}+L}(x). \] 
This definition uses relative distance in $\Grel$.  The limit set of
$\horo{X}{L}$ is the same as the limit set of $X$. This is because if
a sequence $y_n \in \horo{X}{L}$ converges to the Gromov boundary,
then each $y_n$ lies in $B_{\nhat{x_n}+L}(x_n)$ for some $x_n \in
X$. If $y_n$ and $x_n$ limit to distinct points in the Gromov
boundary, then the nearest point projection of $y_n$ to the geodesic
$[x_0, x_n]$ stays a bounded distance from $x_0$, but this implies
that the distance from $y_n$ to $x_n$ is bounded, a contradiction, see
\cite{maher}*{Lemma 3.1} for a more detailed version of this argument.

\begin{theorem} \label{theorem:horoball}
Let $G$ be a weakly relatively hyperbolic group with relative
conjugacy bounds, and let $R$ be a set of elements which are conjugate
to elements of relative length at most $B$. Then for any element $g$, there is a constant
$L$, which only depends on $B$, $\nhat{g}$, and the group constants
$\delta$ and $K$, such that $R \cap Rg$ is contained in an
$L$-horoball neighbourhood of the centralizer of $g$.
\end{theorem}

This theorem shows that the limit set of $R \cap Rg$ is contained in
the limit set of the centralizer of $g$, and hence that the limit set
of $R_k$ is contained in the finite union of limit sets of
centralizers of elements of $G$ with word length at most $k$. In the
case that $R$ is the set of non-\pA elements, the union of the limit
sets of centralizers is dense in the boundary, so the limit set of the
union of the $R_k$, over all $k$, is the entire boundary, which must
be the case, as the limit set of $R$ is the entire boundary.

We start by showing that if $r$ is conjugate to $s$, and the
conjugating element $w$ is chosen to be one of shortest relative
length, then the path in $\Grel$ corresponding to $wsw^{-1}$ is
quasi-geodesic, with quasi-geodesic constants depending only on the
relative length of $s$ and the group constants $\delta$ and $K$, and
independent of the relative length of $r$.

\begin{lemma} \label{lemma:quasigeodesic} 
Let $G$ be a weakly relatively hyperbolic group with relative
conjugacy bounds. Let $r$ be an element of $G$ which is conjugate to
an element $s$, i.e. $r = w s w^{-1}$, for some $w \in G$. If we
choose $w$ to be a conjugating word of shortest relative length, then
the word $wsw^{-1}$ is quasi-geodesic in $\Grel$, with quasi-geodesic
constants which depend only on the relative length of $s$, and the
group constants $\delta$ and $K$. 
\end{lemma}

\begin{proof}
We will use the fact that in a $\delta$-hyperbolic space, a path is a
quasigeodesic if and only if the path lies in a bounded neighbourhood
of a geodesic, and the projection of the path onto the geodesic makes
linear progress along the geodesic. As the path we will consider is a
union of three geodesic segments, one of which has bounded length, it
suffices to show that the path is contained in a bounded neighbourhood
of a geodesic. If the path $wsw^{-1}$ travels far away from a relative
geodesic from $1$ to $r$, then as $s$ has bounded relative length,
there must be a final subsegment of $w$, and an initial subsegment of
$w^{-1}$, which fellow travel. As the path corresponding to $w^{-1}$
is the orientation reverse of $w$, we may choose these initial and
final segments to be inverses of each other. This implies there is a
long final subword of $w$ which conjugates $s$ to an element of
bounded relative length. As the group has relative conjugacy bounds,
we can replace this long subsegment of $w$ with a shorter word,
contradicting our assumption that $w$ was a conjugating element of
shortest relative length. We now write out a detailed version of this
argument.

We will label relative geodesics by their endpoints, so we will write
$[1,w]$ for a particular choice of relative geodesic from $1$ to $w$.
Relative geodesics are in general not unique, but in fact we will not
need to refer to multiple relative geodesics with the same endpoints.
A relative geodesic is a word in the mapping class group, and we will
write $[1,w]^{-1}$ to denote the inverse of this word, which is a
relative geodesic from $1$ to $w^{-1}$. We may also think of the path
$[1,w]^{-1}$ as a translate of $[1,w]$, but with the reverse
orientation. We will also choose relative geodesics $[1,s]$ and
$[1,r]$ from $1$ to $s$ and from $1$ to $r$ respectively.  The group
acts on $\Grel$ on the left by isometries, so $w[1,s]$ is a relative
geodesic from $w$ to $ws$. In particular, the union of the three
relative geodesics $[1,w]$, $w[1,s]$ and $ws[1,w]^{-1}$, is a path
from $1$ to $r$, which corresponds to the concatenation of the words
chosen to represent $w$, $s$ and $w^{-1}$.  This is illustrated below
in Figure \ref{quasi1}.  Note that the segment $ws[1,w]^{-1}$ is equal
to the segment $r[1,w]$, but with the reverse orientation.

\begin{figure}[H]
\begin{center}
\epsfig{file=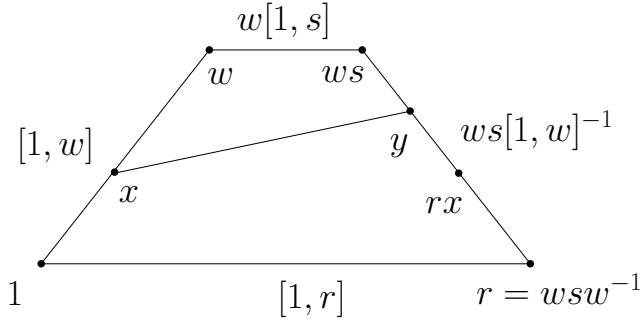, height=120pt}
\end{center}
\caption{The word corresponding to $wsw^{-1}$ gives a path from $1$ to
$r$.} \label{quasi1}
\end{figure}

We may assume that $[1,r]$ and $w[1,s]$ are reasonably far apart. By
thin triangles, the geodesics $[1,w]$ and $ws[1,w^{-1}]$ have long
subsegments which fellow travel. This implies there is a subsegment of
$[1,w]$ which conjugates $s$ to a short word, and as the group has
relative conjugacy bounds, this implies that we may choose $w$ to have
bounded relative length. We now fill in the details of this argument.

Let $D$ be the relative distance from $w[1,s]$ to $[1,r]$. If $D
\leqslant \nhat{s} + 4\delta + 2$, then the entire path $[1,w] \union
w[1,s] \union ws{\geo{w}}^{-1}$ is contained in a $2\nhat{s} + 6\delta
+ 2$ neighbourhood of $\geo{r}$, so assume this is not the case. A
group element $x$ lying on $\geo{w}$ divides the segment into an
initial segment from $1$ to $x$, which we shall denote $\geo{x}$, and
a final segment from $x$ to $w$, which we shall denote
$\path{x}{w}$. Choose $x$ to be the group element of shortest relative
length such that both $\path{x}{w}$ and $r\path{x}{w}$ lie outside a
$2\delta+1$ neighbourhood of $\geo{r}$. Such an element $x$ exists as
the distance from $[1,r]$ to $w[1,s]$ is greater than
$2\delta+1$. Note that the distance from $x$ to $w[1,s]$ is also
greater than $\nhat{s} + 2\delta +1$.

We now show that the distance from $x$ to $rx$ is at most
$\nhat{s}+4\delta$. In a geodesic quadrilateral in a $\delta$-hyperbolic
space, each edge is contained in a $2\delta$ neighbourhood of the
other three, and as $x$ is distance greater than $2\delta$ from either
$w\geo{s}$ or $\geo{r}$, this implies that $x$ lies in a $2\delta$
neighbourhood of the remaining side $ws\geo{w}^{-1}$. Let $y$ be the
closest point on $ws[1,w]^{-1}$ to $x$. By the triangle inequality,
the relative distance from $y$ to $ws$ is at most the relative length
of the path from $y$ to $ws$ via $x$ and $w$.
\begin{align*}
\dhat{y,ws} & \leqslant \dhat{y,x} + \dhat{x,w} + \dhat{w,ws}
\intertext{The distance between $x$ and $y$ is at most $2\delta$, the
  relative distance from $x$ to $w$ is the same as the relative
  distance from $rx$ to $ws$, and the distance between $w$ and $ws$ is
  $\nhat{s}$.}
\dhat{y,ws} & \leqslant 2\delta + \dhat{rx,ws} + \nhat{s} \tag{1}
\intertext{Similarly,
the relative distance from $x$ to $w$ is at most the relative 
length of the path from $x$ to $w$ through $y$ and $ws$.} 
\dhat{x,w} & \leqslant \dhat{x,y} + \dhat{y,ws} + \dhat{ws,w}
\intertext{The relative distance from $x$ to $w$ is the same as the relative
distance from $rx$ to $ws$, the relative distance between $x$ and $y$
is at most $2\delta$, and the length of $\dhat{ws,w}$ is at most
$\nhat{s}$.}
\dhat{rx,ws} & \leqslant 2\delta + \dhat{y,ws} + \nhat{s} \tag{2}
\end{align*}
Together, lines $(1)$ and $(2)$ imply that the difference between the
relative distance of $rx$ from $ws$, and the relative distance of
$y$ from $ws$, is at most $2\delta + \nhat{s}$.  As $rx$ and $y$ lie on a
common relative geodesic through $ws$, this implies that $y$ and $rx$ are
relative distance at most $\nhat{s} + 2\delta$ apart, and hence that $x$ and
$rx$ are relative distance at most $\nhat{s} + 4\delta$ apart.  

The path composed of the three relatively geodesic segments $[x,w]$,
$w[1,s]$ and $r[w,x]$ gives a word in the mapping class corresponding
to $(x^{-1}w)s(x^{-1}w)^{-1}$. The relative length of this group
element is the distance between $x$ and $rx$, which is at most
$\nhat{s}+4\delta$, so $x^{-1}w$ conjugates $s$ to a word of length at
most $\nhat{s}+4\delta$. As we have assumed that the group $G$ has
relative conjugacy bounds, we may choose the relative length of the
conjugating element $x^{-1}w$ to be at most $K(2\nhat{s} + 4\delta)$,
where $K$ is the relative conjugacy bound constant.  This implies that
distance $D$ between $w[1,s]$ and $[1,r]$ is at most $K(\nhat{s} +
2\delta) + 2\delta+1$.  Hence the union of the segments $[1,w]$,
$w[1,s]$ and $ws[1,w]^{-1}$ lies in an $L_1$-neighbourhood of $[1,r]$,
where $L_1 = \max \{ K(\nhat{s} + 2\delta) + 2\delta+1 +\nhat{s},
2\nhat{s} + 4\delta + 2 \} \le 2K(\nhat{s} + 2\delta) + 2$. The
constant $L_1$ depends only on $\nhat{s}$, and the group constants
$\delta$ and $K$.
\end{proof}

We now observe that if $r$ is conjugate to $s$ by a word $w$ of shortest
relative length, then the relative length of $w$ is roughly half the
relative length of $r$.

\begin{proposition} \label{proposition:halfway} Let $G$ be a
  relatively hyperbolic group.
Let $r$ in $G$ be conjugate to $s$ by $w$, such that $w$ and $ws$ are contained in
an $L_1$-neighbourhood of a relative geodesic $[1,r]$. Then \[
\frac{1}{2}(\nhat{r} - \nhat{s}) \leqslant \nhat{w} \leqslant
\frac{1}{2}(\nhat{r} + \nhat{s}) + 2L_1. \]
\end{proposition}

\begin{proof}
By the triangle inequality $\nhat{r} \leqslant 2\nhat{w}+\nhat{s}$,
giving the left hand inequality. Let $a$
be the closest point on $[1,r]$ to $w$, and let $b$ be the closest
point on $[1,r]$ to $ws$. Relative geodesics $[w,a]$ and $[ws,b]$ have length
at most $L_1$. Let $[1,a]$ be the initial segment of $[1,r]$ from $1$
to $a$, and let $[b,r]$ be the terminal segment of $[1,r]$ from $b$ to
$r$. The overlap of $[1,a]$ and $[b,r]$ is at most the relative 
length of $[a,b]$, which is at most $2L_1 + \nhat{s}$, giving $\nhat{a} +
\dhat{b,r} \leqslant \nhat{r} + 2L_1 + \nhat{s}$. Using the triangle inequality, $\nhat{w} 
\leqslant \nhat{a} + L_1$ and $\nhat{w} \leqslant \dhat{b,r} +
L_1$, which
implies $2\nhat{w} \leqslant \nhat{r} + 4L_1 + \nhat{s}$. This gives the
right hand inequality.
\end{proof}

We can now complete the proof of Theorem \ref{theorem:horoball}.

%

\begin{proof}
Suppose $r' = rg$, and $r$ and $r'$ are conjugate to relatively short
words $s$ and $s'$ respectively. By Lemma \ref{lemma:quasigeodesic},
we may choose shortest conjugating words $w$ and $w'$ so that the
paths $wsw^{-1}$ and $w's'w'^{-1}$ are quasigeodesic, with endpoints
close together. This implies that $w$ and $w'$ are close together, and
as the relative lengths of $s$ and $s'$ are bounded, this implies that
$w^{-1}$ conjugates $g$ to another short word $w^{-1}gw$, and so $w$ is close to the
centralizer of $g$. This in turn implies that $r$ lies in a horoball
neighbourhood of the centralizer $C(g)$. We now give a detailed
version of this argument.

It suffices to prove the result for elements $r \in R \cap Rg$ with
$\nhat{r}$ sufficiently large. If $r$ is in $R \cap Rg$, then there is
an $r' \in R$ such that $r^{-1}r' = g$. Let $[1,r]$ be a relative
geodesic from $1$ to $r$, and let $[1,r']$ be a relative geodesic from
$1$ to $r'$. By thin triangles, the two relative geodesics $[1,r]$ and
$[1,r']$ are relative distance at most $\delta$ apart for initial
segments of relative length at least $\nhat{r} - \nhat{g} - \delta$.
Later on in this argument we will need this to be be at least
$\frac{1}{2}(\nhat{r} + 4L_1 + B) + L_1 + \delta$, where $L_1$ is the
constant from Lemma \ref{lemma:quasigeodesic} above, so we will assume
that $\nhat{r} > 6L_1 + B + 4\delta + 2\nhat{g}$. In particular, this
implies that $L$ must be at least as large as this value.

Let $w$ be a word of shortest relative length conjugating $r$ to an
element $s$ of relative length at most $B$, and let $w'$ be a word of
shortest relative length conjugating $r'$ to a word $s'$ of relative
length at most $B$. Choose relative geodesics $[1,w]$ and $[1,s]$ from
$1$ to $w$ and $s$ respectively. Then the union of the three
relatively geodesic segments $[1,w] \cup w[1,s] \cup ws[1,w]^{-1}$ is
a path from $1$ to $r = wsw^{-1}$. Similarly, choose relative
geodesics $[1,w']$ and $[1,s']$ from $1$ to $w'$ and $s'$
respectively. Then the union of the three relatively geodesic segments
$[1,w'] \cup w'[1,s'] \cup w's'[1,w']^{-1}$ is a path from $1$ to $r'
= w's'w'^{-1}$. This is illustrated below in Figure \ref{quasi2}.

\begin{figure}[H]
\begin{center}
\epsfig{file=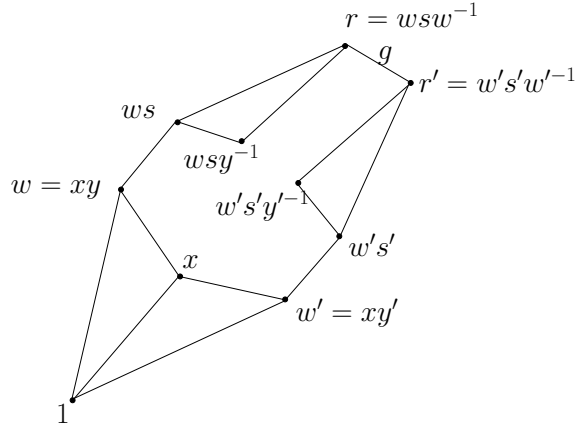, height=160pt}
\end{center}
\caption{The group elements $r$ and $r'$ are conjugate to $s$ and $s'$
respectively.} \label{quasi2}
\end{figure}

We now show that the two group elements $w$ and $w'$ are close
together. By Proposition \ref{proposition:halfway}, the relative
length of $w$ is roughly half that of $r$, i.e.
$\frac{1}{2}(\nhat{r}-B) \leqslant \nhat{w} \leqslant
\frac{1}{2}(\nhat{r}+4L_1+B)$, and there are similar bounds for $w'$.
By Lemma \ref{lemma:quasigeodesic}, the group element $w$ lies in an
$L_1$-neighbourhood of the geodesic $[1,r]$.  Let $a$ be the closest
point on $[1,r]$ to $w$, then the relative distance between $w$ and
$a$ is at most $L_1$, so $\frac{1}{2}(\nhat{r}-B) - L_1 \leqslant
\nhat{a} \leqslant \frac{1}{2}(\nhat{r}+4L_1+B) + L_1$. Again, by
Lemma \ref{lemma:quasigeodesic}, the group element $w'$ lies in an
$L_1$-neighbourhood of $[1,r']$, which, as we have assumed that
$\nhat{r}$ is sufficiently large, has an initial segment of at least
this length contained in a $\delta$ neighbourhood of $[1,r]$. So if
$a'$ is the closest point on $[1,r]$ to $w'$, then the distance from
$a'$ to $w'$ is at most $L_1 + \delta$. This implies that
$\frac{1}{2}(\nhat{r}-B) - L_1 - \delta \leqslant \nhat{a'} \leqslant
\frac{1}{2}(\nhat{r}+4L_1+B) + L_1 + \delta$. Therefore, as $a$ and
$a'$ lie on a common relative geodesic through $1$, the relative
distance between $a$ and $a'$ is at most $4L_1 + B + \delta$, so the
relative distance between $w$ and $w'$ is at most $6L_1 + 2B +
2\delta$.

As $w$ and $w'$ are close together, we may choose words to represent
them that have a large common initial segment consisting of a relative
geodesic $[1,x]$. The relative distance between $w$ and $w'$ is at
most $6L_1 + 2B + 2\delta$, so we may choose a group element $x$ such
that $w = xy$ and $w'=xy'$, for $y$ and $y'$ of relative length at
most $3L_1 + B + \delta$. The path from $ws$ to $r$ gives a word
corresponding to $w^{-1}$; in fact this path is precisely
$ws[1,w]^{-1}$. As $w^{-1} = y^{-1}x^{-1}$, we may also travel from
$ws$ to $r$ using the path $ws[1,y]^{-1} \cup wsy^{-1}[1,x]^{-1}$.
Similarly, the path from $w's'$ to $r'$ is $w's'[1,w']^{-1}$, so we
may also travel from $w's'$ to $r'$ along the path $w's'[1,y']^{-1}
\cup w's'y'^{-1}[1,x]^{-1}$.  This is illustrated in Figure
\ref{quasi2} above.

We now show that the group element $x$ conjugates $g$ to an element of
bounded relative length. The two group elements $wsy^{-1}$ and
$w's'y'^{-1}$ are endpoints of a path consisting of three relatively
geodesic segments passing through $r$ and $r'$, namely
$wsy^{-1}[1,x]^{-1} \cup r[1,g] \cup r'[1,x]$.  The group element
corresponding to this path is $x^{-1}gx$, which is a conjugate of $g$,
and we now show that the relative length of this element is bounded,
i.e. we show that the endpoints of the path are a bounded relative
distance apart. The relative distance between $w$ and $w'$ is at most
$6L_1 + 2B + 2\delta$. The group elements $s$ and $s'$ have relative
length at most $B$, so the distance between $ws$ and $w's'$ is at most
$6L_1 + 4B + 2\delta$. The elements $y$ and $y'$ have length at most
$3L_1 + B + \delta$, so the relative distance between $wsy^{-1}$ and
$ws'y'^{-1}$ is at most $12L_1 + 6B + 4\delta$.

Finally, we show that $w$ lies close to the centralizer of $g$, which
will imply that $r$ lives in a horoball neighbourhood of $C(g)$.  As
we have assumed that the group $G$ has relative conjugacy bounds,
there is a group element $v$ of length at most $K(\nhat{g}+12L_1 + 6B
+ 4\delta)$ which conjugates $x^{-1}gx$ to $g$, where $K$ is the
relative conjugacy bound constant. This implies that $x^{-1}gx =
vgv^{-1}$. Multiplying on the right by $v$ and on the left by $x$, we
obtain $gxv = xvg$, so $xv$ lies in the centralizer of $g$. Therefore
$x$ is a distance at most $K(\nhat{g}+12L_1 + 6B + 4\delta)$ from
$C(g)$, and so $w$ is a distance at most $L_2 = K(\nhat{g}+12L_1 + 6B
+ 4\delta) + 3L_1 + B + \delta$ from $C(g)$. Let $b$ be the closest
point in $C(g)$ to $w$.  The triangle inequality implies that
$\nhat{w} \le \nhat{b} + L_2$.  The distance from $wsw^{-1}$ to $C(g)$
is at most $\nhat{w} + B + L_2$, which is at most $\nhat{b} + 2L_2
+B$. Therefore $wsw^{-1}$ lies in an $L$-horoball neighbourhood of
$C(g)$, where $L = 2L_2 + B$.  The constant $L$ depends only on $B$,
$\nhat{g}$, and the group constants $\delta$ and $K$.
\end{proof}

This completes the proof of Theorem \ref{theorem:horoball}.

\section{Random walks} \label{section:random}

We start by recalling some basic definitions about random walks on
groups, see for example Woess \cite{woess}. We then show that a random
walk on the mapping class group converges to a projective measured
lamination with probability one, and observe that the periodic and
reducible elements are all conjugate to elements of bounded relative
length, for some bound which depends on the surface.  Then in Section
\ref{section:measure} we show that the limit set of the centralizer of
an element of the mapping class group has harmonic measure zero, as
long as the limit set has infinitely many images under the group $H$
generated by the support of the random walk. Let $R$ be a set of
elements of $G$ which are conjugate to elements of bounded relative
length, and let $R_k$ be the set of $k$-dense elements of $R$. We
showed in Section \ref{section:coarse} that the limit set of $R_k$ is
contained in a finite union of limit sets of non-trivial centralizers,
so this implies that the limit set of $R_k$ has harmonic measure zero,
under the assumption on the images of limit sets of centralizers.
Finally, in Section \ref{section:asymptotic} we argue that this means
a sample path travels through regions where the elements of $R$ become
further and further apart, so the probability that the random walk is
in $R$ tends to zero asymptotically. We finally complete the proof by
showing how to reduce the general case to the case in which all limit
sets of centralizers of non-trivial elements have infinitely many
images under $H$.

Let $G$ be a group, and let $\mu$ be a probability distribution on
$G$, which we shall call the \emph{step distribution}. We may use the
probability distribution $\mu$ to generate a Markov chain, or
\emph{random walk} on $G$, with transition probabilities $p(x,y) =
\mu(x^{-1}y)$, and we shall always assume that we start at time zero
at the identity.  We shall write $p^{(n)}(x,y)$ for the probability
that you go from $x$ to $y$ in $n$ steps.  The \emph{step space} for
the random walk is the infinite product $(G, \mu)^{\Z_+}$, i.e. the
\emph{steps} or \emph{increments} of the random walk are a sequence of
independent identically $\mu$-distributed random variables. An element
of the step space $(s_1, s_2, s_3, \ldots )$ determines a path in $G$,
where the location $w_n$ of the path at time $n$ is given by the
product of the first $n$ steps, i.e. $w_n = s_1 s_2 \dots s_n$. This
gives a map $\rho : G^{\Z_+} \to G^{\Z_+}$ from steps to paths, and we
will write $(G^{\Z_+},\P_\mu)$ for the \emph{path space}, where
$\P_\mu(X)$ is the measure of $\rho^{-1}(X)$ in the step space. We
will often just write $\P$ for $\P_\mu$ if it is clear from context
which probability distribution determines the random walk. We will
call an element $\w = (w_1, w_2, \ldots)$ of the path space a
\emph{sample path}. The distribution of random walks at time $n$ is
given by the $n$-fold convolution of $\mu$, which we shall write as
$\mu^{(n)}$. We shall always require that the group generated by the
support of $\mu$ is non-elementary. We do not assume that the
probability distribution $\mu$ is symmetric, so the group generated by
the support of $\mu$ may be strictly larger than the semi-group
generated by the support of $\mu$.

A homomorphism $\phi \colon G \to H$ takes random walks on $G$ to
random walks on $H$. More precisely, the homomorphism $\phi$
determines a map, which we shall also call $\phi$, from $G^{\Z_+}$ to
$H^{\Z_+}$, which sends a sequence $(s_1, s_2, s_3, \ldots)$ to the
sequence $(\phi(s_1), \phi(s_2), \phi(s_3), \ldots)$. Let $\phi^* \mu$
be the probability distribution on $H$ defined by $\phi^* \mu(h) =
\mu(\phi^{-1}(h))$, then $\phi$ gives a measure preserving map from
$(G, \mu)^{\Z_+}$ to $(H, \phi^* \mu)^{\Z_+}$, such that the following
diagram commutes,
\[ \begin{CD}
(G,\mu)^{\Z_+} @>\phi>> (H, \phi^* \mu)^{\Z_+} \\
@VV{\rho}V                @VV{\rho}V \\
(G^{\Z_+}, \P_\mu) @>\phi>> (H^{\Z_+}, \P_{\phi^* \mu}) \\
\end{CD} \] and so the lower map is also measure preserving. Therefore
the images of sample paths in $H$ under $\phi$ are distributed
according to a random walk on $H$ determined by the probability
distribution $\phi^* \mu$.

An example of a random walk is the nearest neighbour random
walk on a Cayley graph $\G$ for $G$. This is the random walk
determined by a probability distribution $\mu$ which gives equal
weight to every generator and its inverse, and is zero on all the
other elements of the group. In this case the support of the random
walk is the whole group.

We now show that a sample path $\w$ gives a sequence of points in
$\Grel$, which converges to a lamination $\l(\w)$ in the Gromov
boundary of the complex of curves, almost surely. Given a subset $X$
of $\lmin$, the measure of $X$ with respect to $\nu$ is defined to be
the proportion of sample paths which converge to points contained in
$X$. This measure $\nu$ is $\mu$-stationary, i.e., for any subset $X$
of $\lmin$, \[ \nu(X) = \sum_{g \in G} \mu(g) \nu(g^{-1}X). \] We will
also refer to $\nu$ as a harmonic measure on $\lmin$.

\begin{theorem} 
\label{theorem:convergence}
Consider a random walk on the mapping class group of a non-sporadic
orientable surface of finite type, determined by a probability
distribution $\mu$ such that the group generated by the support of
$\mu$ is non-elementary. Then the sequence of points $\{ w_n \}$
determined by a sample path $\w$ converge to the boundary of the
relative space $\Grel$ with probability one, and the distribution of
points on the boundary is given by a unique $\mu$-stationary measure
$\nu$ on $\lmin$, the Gromov boundary of $\Grel$.
\end{theorem}

Theorem \ref{theorem:convergence} is an immediate consequence of
results of Kaimanovich and Masur \cite{km}, and Klarreich
\cite{klarreich}, which we now describe.  Kaimanovich and Masur
consider the action of the mapping class group on Teichm\"uller space,
and show that given a basepoint $x_0 \in \T(\S)$, then for almost all
sample paths, the sequence of images of the basepoint $\{ w_n x_0\}$
converges to a point in $\PMF$, the Thurston compactification of
Teichm\"uller space.

\begin{theorem}{\cite{km}*{Theorem 2.3.4}} \label{theorem:km} 
If $\mu$ is a probability measure on the mapping class group $G$ such
that the group generated by its support is non-elementary, then there
exists a unique $\mu$-stationary probability measure $\nu$ on the
space $\PMF$, which is purely non-atomic, and concentrated on the
subset $\UE \subset \PMF$ of uniquely ergodic foliations.
For any $x \in \T(\S)$ and almost every sample path $\omega = \{ w_n
\}$ of the random walk determined by $(G,\mu)$, the sequence $w_nx$
converges in $\PMF$ to a limit $F(\omega) \in \UE$, and the
distribution of the limits $F(\omega)$ is given by $\nu$.
\end{theorem}

Kaimanovich and Masur \cite{km} only state the result for closed
surfaces, but Farb and Masur \cite{fm} point out that the proof works
for surfaces with punctures.

The following result of Klarreich \cite{klarreich} describes the
Gromov boundary of the complex of curves, and will enable us to relate
convergence in $\PMF$ to convergence in the Gromov boundary of the
complex of curves.

\begin{theorem}{\cite{klarreich}*{Theorem 1.2}}
The inclusion map from $\T(\S)$ to $\T_{el}(\S)$ extends continuously to
the portion $\fmin(\S)$ of $\PMF(\S)$ consisting of minimal
foliations, to give a map $\pi:\fmin(\S) \to
\d\T_{el}(\S)$. The map $\pi$ is surjective, and $\pi(F) =
\pi(G)$ if and only if $F$ and $G$ are topologically
equivalent. Moreover any sequence $\{x_n\}$ in $\T(\S)$ that converges
to a point in $\PMF(\S) \setminus \fmin(\S)$ cannot accumulate in the
electrified space onto any portion of $\d\T_{el}(\S)$.  
\end{theorem}

A sample path $\w$ gives rise to a sequence of points $\{ w_n x_0 \}$
in Teichm\"uller space which converge to a uniquely ergodic foliation
$F(\w)$, almost surely. Uniquely ergodic foliations are minimal, so
the image of the sequence of points in $\T_{el}$ converges to the same
foliation $F(\w)$. As $\T_{el}$ is quasi-isometric to $\C(\S)$, the
sequence $\{ w_n x_0 \}$ in $\T_{el}$ gives rise to a sequence in
$\C(\S)$, which lies a bounded distance from the images of $w_n$ in
$\Grel$. Therefore the sequence $w_n$ in $\Grel$ converges almost
surely to a point in the Gromov boundary corresponding to a uniquely
ergodic foliation $F(\w)$. This enables us to define a harmonic
measure on the Gromov boundary $\fmin$, however, as the map from
$\PMF$ to $\fmin$ is a bijection when restricted to the uniquely
ergodic foliations, this harmonic measure is the same as the pullback
of $\nu$ under the map $\PMF \to \fmin$. We will abuse notation and
write $\nu$ for the harmonic measure on either $\PMF$ or $\fmin$. As
$\PMF$ is essentially the same as $\PML$, and $\fmin$ is essentially
the same as $\lmin$, we will also write $\nu$ for harmonic measure on
these spaces as well.

The rest of this section is devoted to the proof of Theorem
\ref{theorem:rw}. In fact, we will prove the following more general
result.

\begin{theorem}
\label{theorem:main} 
Consider a random walk on the mapping class group of an orientable
non-sporadic surface of finite type, determined by a probability
distribution $\mu$, whose support generates a non-elementary subgroup. 
Let $R$ be a subset of $G$ with the property that every element of $R$
is conjugate to an element of relative length at most $B$, for some
constant $B$. Then the
probability that a random walk of length $n$ lies in $R$ tends to zero
as $n$ tend to infinity.
\end{theorem}

The following observation shows that the set of periodic and reducible
elements of the mapping class group form a set of elements which are
conjugate to elements of bounded relative length, so Theorem
\ref{theorem:main} implies Theorem \ref{theorem:rw}.

\begin{lemma} \label{lemma:bounded length}
Every reducible or periodic element of the mapping class group of a
surface which is not a sphere with three or fewer punctures is
conjugate to an element of bounded relative length, where the bound
only depends on the surface $\S$.
\end{lemma}

\begin{proof}
If $g$ is reducible, then $g$ preserves a collection of disjoint
simple closed curves. We can conjugate $g$ so that one of these curves
$y$ is distance $1$ from the basepoint $x_0$ in $\C(\S)$. As $g(y)$ is
disjoint from $y$, and hence distance at most $1$ from $y$, the image
of the basepoint $x_0$ under $g$ is distance at most $3$ from $x_0$,
so $g$ is conjugate to an element of relative length at most $3\Kcc$,
where $\Kcc$ is the quasi-isometry constant between the relative
metric and the complex of curves.

There are only finitely many conjugacy classes of periodic elements in
the mapping class group of a given surface, so every periodic element
is conjugate to a periodic element of bounded relative length, for
some bound that depends on the surface.
\end{proof}

We now prove Theorem \ref{theorem:main}. In Section
\ref{section:measure} we show that the harmonic measure of the limit
set of a centralizer is zero, assuming that there are infinitely many
images of the limit set under the group generated by the support of
the random walk. Finally in Section \ref{section:asymptotic}, we show
that this implies that the asymptotic probability of the random walk
being in $R$ tends to zero.

\subsection{Centralizers have harmonic measure zero} \label{section:measure}

In this section we show that if there are infinitely many images of the
limit set of a centralizer $C(g)$ under the group generated by the
support of the random walk, then $\overline{C(g)}$ has harmonic
measure zero.

\begin{lemma} 
\label{lemma:small} 
Consider a random walk on the mapping class group $G$ of an orientable
non-sporadic surface, determined by a probability distribution $\mu$,
whose support generates a non-elementary subgroup $H$. Then if the
limit set of a centralizer $\overline{C(g)}$ has infinitely many
images under $H$, then  $\nu(\overline{C(g)}) =
0$, where $\nu$ is the harmonic measure
on $\lmin$ determined by $\mu$.
\end{lemma}

We do not assume that the probability distribution $\mu$ is symmetric,
so the group $H$ generated by the support of $\mu$ may be larger than
the semi-group generated by the support of $\mu$, which we shall
denote by $H^+$. We shall write $H^-$ for the semi-group generated by
the inverses of elements in the support of $\mu$.

\begin{proof}
If $g$ is not periodic, then $\overline{C(g)}$ contains
at most two points in the Gromov boundary $\lmin$, by Proposition
\ref{prop:small limit set}, so $\overline{C(g)}$ has harmonic measure
zero, as $\nu$ is non-atomic. We now consider the case in which $g$ is
periodic, and it will then be convenient to consider arbitrary finite
subgroups of the mapping class group, rather than just cyclic
subgroups.  We will argue that we may choose a finite subgroup $F$
such that the intersections of infinitely many distinct images of its
limit sets have measure zero.  The $\mu$-stationarity of the measure
$\nu$ will then imply that all of these limit sets have $\nu$-measure
zero.

Let $F$ be a maximal finite subgroup such that $\nu(\overline{C(F)}) >
0$. We now show that this implies that if any two distinct images of
$\overline{C(F)}$ intersect, their intersection has measure zero with
respect to $\nu$. The intersection of two images of $C(F)$ is the
centralizer of the subgroup $F'$ generated by the two conjugates of
$F$, which is strictly larger than $F$ as the images are distinct. If
the subgroup $F'$ is finite, then it has a limit set of measure zero,
as we have assumed that $F$ is a maximal finite subgroup with
$\nu(\overline{C(F)}) > 0$. If the subgroup $F'$ is infinite, then by
Proposition \ref{prop:small limit set}, the limit set of the
centralizer of $F'$ consists of at most two points, and so has measure
zero, as $\nu$ is non-atomic.

Let $s$ be the supremum of the harmonic measure of images of
$\overline{C(F)}$ under the action of $H^-$, i.e.  $s = \sup_{h \in
  H^-} \nu(h\overline{C(F)})$. The semi-group $H^-$ generates $H$, so
by Proposition \ref{prop:infinite images}, there are infinitely many
distinct images of $\overline{C(F)}$ under the action of $H^-$, and
furthermore, these images intersect in smaller centralizers, which
have $\nu$-measure zero, by our assumption on the maximality of $F$.  First
suppose that the supremum is achieved.  If $\nu(h\overline{C(F)})$ is
equal to the supremum for some $h$, then as $\nu$ is $\mu$-stationary,
this means that
\[ s = \nu(h\overline{C(F)}) = \sum_{g \in H^+} \mu(g)
\nu(g^{-1}h\overline{C(F)}). \] 
We may take the sum over $H^+$ as $\mu(g) = 0$ for $g \not \in H^+$.
As $g^{-1}$ is in $H^-$, this implies that
$\nu(g^{-1}h\overline{C(F)}) \le s$, but as $\mu$ is a probability
distribution with total mass one, if any $\nu(g^{-1}h\overline{C(F)})$
is strictly less than $s$ for $\mu(g) \not = 0$, then the right hand
sum is strictly less than $s$, a contradiction. The harmonic measure
$\nu$ is $\mu^{(n)}$-stationary, for all $n$, so this implies that
$\nu(gh\overline{C(F)}) = s$ for all elements $g$ in the semi-group
$H^-$, which contradicts the fact that $\nu$ has finite total mass.
If the supremum is not achieved, then let $h_i$ be a sequence of
elements in $H^-$ such that $\nu(h_i\overline{C(F)})$ tends to the
supremum $s$, and we may assume we have chosen a sequence in which the
$h_i \overline{C(F)}$ are all distinct. This gives infinitely many
distinct sets each with measure bounded away from zero, and pairwise
intersections of measure zero, which contradicts the fact that $\nu$
is a measure with finite total mass.

This completes the proof of Lemma \ref{lemma:small}.
\end{proof}

\subsection{Asymptotic probabilities} \label{section:asymptotic}

A random walk is \emph{recurrent}
on a subset $X$ of $G$, if a sample path hits $X$ infinitely often
with probability one.  We say a random walk is \emph{transient} on a
subset $X$ of $G$, if a sample path hits $X$ infinitely often with
probability zero.

We now observe that if the harmonic measure of the limit set of a set
$X$ is zero, then the set $X$ is transient. In particular the
probability that a random walk of length $n$ is in $X$ tends to zero
as $n$ tend to infinity.

\begin{lemma}
  Let $X$ be a subset of the mapping class group of a non-sporadic
  surface. If the harmonic measure of the closure of $X$ is zero, then
  the random walk is transient on $X$.
\end{lemma}

\begin{proof}
If a sample path converges to a lamination, and is also recurrent on
$X$, i.e.  hits $X$ infinitely often, then the limiting lamination
must lie in the closure $\overline X$. As sample paths converge with
probability one, the probability that a sample path is recurrent on $X$ is
bounded above by the harmonic measure of $\overline X$. In
particular. if the harmonic measure of $\overline X$ is zero, then the
random walk is transient on $X$, i.e. a sample path hits $X$ finitely
many times with probability one.
\end{proof}

A sample path which never hits $X$ may still converge to a lamination
lying in $\overline X$, so if the harmonic measure of $\overline X$ is
greater than zero, this does not imply that sample paths are recurrent
on $X$ with probability greater than zero.

We now prove Theorem \ref{theorem:main}, under the following assumption.

\begin{description}
\item[$(*)$] The limit set of the centralizer of every non-trivial element of $G$ has
infinitely many images under the group $H$ generated by the support of
the random walk. 
\end{description}

Recall that $R_k$ is the set of $k$-dense elements of $R$, so any two
points in $R \setminus R_k$ are distance at least $k$ apart. We shall
say $R \setminus R_k$ is a \emph{$k$-separated} set. Furthermore, we
have shown that the limit set of $R_k$ is contained in the limit set
of a finite union of centralizers. As we have assumed that every
$\overline{C(g)}$ has infinitely many images under $H$, Lemma
\ref{lemma:small} implies that each $\overline{C(g)}$ has harmonic
measure zero, and so the limit set of $R_k$ also has harmonic measure
zero.  We want to show that the probability that you lie in $R$ in a
low density region is small. The basic idea is that outside of $R_k$,
the distance between elements of $R$ is at least $k$, so for a sample
path outside $R_k$ there is some upper bound, depending on $k$, for
how often the sample path hits elements of $R$, and furthermore, this
upper bound tends to zero as $k$ tends to infinity.

\begin{figure}[H] 
\begin{center}
\epsfig{file=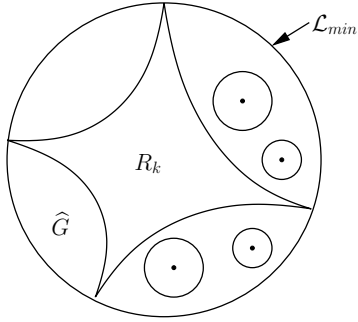, height=120pt}
\end{center}
\caption{$R_k$ has harmonic measure zero.}\label{picture30}
\end{figure}

\begin{lemma} \label{lemma:ref}
Let $X$ be a $k$-separated set in $G$. Then 
\[ \mun(X) \le \max \{ \mu^{(m)}(g) \mid g \in B_{k/2} \} +
\mu^{(m)}(G \setminus B_{k/2}) \]
for all $m < n$, where $B_{k/2}$ is the ball of radius $k/2$ about
the origin in the word metric.
\end{lemma}

\begin{proof}
The measure $\mun$ is the convolution of $\mu^{(n-m)}$ and
$\mu^{(m)}$, i.e.
\[ \mun(X) = \sum_{g \in G} \mu^{(n-m)}(g) \mu^{(m)}(g^{-1}X), \]
for any $m < n$. Any translate $gX$ of $X$ is also $k$-separated, so at most one point
of $gX$ intersects the ball of radius $k/2$ about the origin, therefore
\[\mu^{(m)}(gX) \le \max\{ \mu^{(m)}(g) \mid g \in B_{k/2} \}  + \mu^{(m)}(G
\setminus B_{k/2}) \] for all $m$, and any $g \in G$. This implies
$\mun(X) \le \max\{ \mu^{(m)}(g) \mid g \in B_{k/2} \}  + \mu^{(m)}(G
\setminus B_{k/2})$, as required.
\end{proof}

Let $s_m = \sup \{ \mu^{(m)}(g) \mid g \in G \}$, which is an upper
bound for $\max \{ \mu^{(m)}(g) \mid g \in B_{k/2} \}$. We now show
that $s_m$ tends to zero as $m$ tends to infinity. If $s_m$ does not
tend to zero, then there is a sequence $g_m$ with $\mu^{(m)}(g_m) \ge
\e > 0$. Therefore there is a set of sample paths of positive measure
(in fact of measure at least $\e$) which are recurrent on any infinite
subset of $\{ g_m \}$. Choose a basepoint $x_0$ in Teichm\"uller space
$\T(\S)$, and consider the sequence of images $g_m(x_0)$. As sample
paths converge to the boundary almost surely, the random walk is
transient on bounded sets, so the distance in the Teichm\"uller metric
between $x_0$ and $g_m(x_0)$ tends to infinity. Therefore we may pass
to a subsequence such that $g_m(x_0)$ converges to a foliation $F$ in
$\PML$. As $\nu$ is the weak-$*$ limit of $\mu^{(n)}$ on $\T(\S) \cup
\PML$, this implies that the harmonic measure of $F$ is strictly
larger than zero, contradicting the fact that the harmonic measure
$\nu$ is non-atomic.  This shows that $s_m$ tends to zero as $m$ tends
to infinity.

To summarize, we have shown
\begin{align*}
\P(w_n \in R) & = \P(w_n \in R_k) + \P(w_n \in R \setminus R_k) \\
& \leqslant \mun(R_k) + s_m  + \mu^{(m)}(G \setminus B_{k/2})
\intertext{%
for all $k$, and any $m < n$. The limit sets of centralizers have
measure zero, so the harmonic measure of $R_k$ is zero, so $\mun(R_k)$
tends to zero as $n$ tends to infinity, for every $k$ and
$m$. Therefore
}
\lim_{n \to \infty} \P(w_n \in R) & \leqslant s_m + \mu^{(m)}(G \setminus B_{k/2}),
\end{align*}
for all $k$ and $m$. We have shown that $s_{m}$ tends to zero as $m$
tends to infinity, and we can choose a sequence $k_m$, which also
tends to infinity as $m$ tends to infinity, such that $\mu^{(m)}(G
\setminus B_{k_m/2})$ tends to zero. This implies that the probability
that $w_n$ lies in $R$ tends to zero as $n$ tends to infinity, as
required. This completes the proof of Theorem \ref{theorem:main},
assuming property $(*)$ above. 

Finally, we prove Theorem \ref{theorem:main} when there are limit sets
of centralizers with finitely many images under the subgroup generated
by the support of the random walk. We start by showing that the
surface $\S$ covers a surface $\O$ such that there is a homomorphism
$\phi$ from $H$ to $G_\O$ with finite kernel, such that the image of
$H$ in $G_\O$ satisfies property $(*)$, i.e. every non-trivial finite
subgroup $\overline{C(F)}$ has infinitely many images under
$\phi(H)$. Suppose there is a finite subgroup $F_1$ in $G$ such that
$\overline{C(F_1)}$ has only finitely many images under $H$. Then by
Proposition \ref{prop:infinite images}, the subgroup $H$ is contained
in the normalizer $N(F_1)$, and so there is a map $\phi_1$ from $H$ to
$G_{\O_1}$ with finite kernel, where $\O_1$ is the quotient surface
$\S / F_1$. We may repeat this process, i.e. if there is a
non-trivial finite subgroup $F_2$ in $G_{\O_1}$ such that
$\overline{C(F_2)}$ has only finitely many images under $\phi_1(H)$,
then $\phi_1(H)$ is contained in $N(F_2)$, and so we obtain a
homomorphism $\phi_2 \circ \phi_1$ from $H$ to $G_{\O_2}$ with finite
kernel, where $\O_2 = \O_1 / F_2$. However, the area of the hyperbolic
surface is divided by the degree of the covering at each stage, and so
this process terminates after finitely many steps, as there is a lower
bound on the area of a hyperbolic orbifold. This gives the desired
homomorphism from $H$ to $G_\O$ with finite kernel, for some surface
$\O$ covered by $\S$.

The homomorphism $\phi$ has finite kernel, so $\phi(H)$ is
non-elementary subgroup of $G_\O$. The relative metric on $G_\O$ is
quasi-isometric to the complex of curves $\C(\O)$, and the map on
curve complexes from $\C(\O)$ to $\C(\S)$ induced by the covering is a
quasi-isometric embedding, by Rafi and Schleimer \cite{rs}, so the
relative metric on $H$ is also quasi-isometric to the relative metric
on $\phi(H)$. Therefore the image of $\phi(R \cap H)$ is contained in
a set $R' \subset G_\O$ of elements which are conjugate to elements of
bounded relative length, possibly for a different bound on the
relative length.  The homomorphism $\phi$ gives rise to an induced
random walk on $G_\O$ determined by the probability distribution
$\phi^* \mu(h) = \mu(\phi^{-1}(h))$.  Therefore, as the induced random
walk on $G_\O$ satisfies property $(*)$, the probability that the
image of the random walk lies in $R'$ tends to zero, and so the
probability that the original random walk lies in $R$ also tends to
zero. This completes the proof of Theorem \ref{theorem:main}.



\end{document}